\documentclass[12pt, leqno]{amsart}
\usepackage{amssymb,amsmath, amscd}

=\oddsidemargin

\font\teneufm=eufm10
\font\seveneufm=eufm7
\font\fiveeufm=eufm5
\newfam\eufmfam
\textfont\eufmfam=\teneufm
\scriptfont\eufmfam=\seveneufm 
\scriptscriptfont\eufmfam=\fiveeufm

\let\goth\mathfrak

\def\gg{\goth g}
\def\gs{\goth s}
\def\gh{\goth h}

\def\gl{\goth l}

\def\gp{\goth p}

\def\gn{\goth n}

\def\gq{\goth q}

\def\go{\goth o}
\def\ga{\goth a}

\def\beq{\begin{equation}}
\def\eeq{\end{equation}}
\def\bea{\begin{eqnarray}}
\def\eea{\end{eqnarray}}
\def\beas{\begin{eqnarray*}}
\def\eeas{\end{eqnarray*}}

\def\cplus{\hbox{$\supset${\raise1.05pt\hbox{\kern -0.55em
${\scriptscriptstyle +}$}}\ }}

\DeclareMathOperator{\Span}{Span}
\DeclareMathOperator{\tr}{tr}

\DeclareMathOperator{\Id}{Id}

\def \wt {\widetilde}

\newtheorem{theorem}[equation]{Theorem}

\newtheorem{lemma}[equation]{Lemma}

\newtheorem{corollary}[equation]{Corollary}

\newtheorem{proposition}[equation]{Proposition}

\newtheorem{definition}[equation]{Definition}

\newtheorem{remark}[equation]{Remark}

\def\Z{\mathbb Z}
\def\C{\mathbb C}
\def\R{\mathbb R}

\begin{document}

\title[Cominuscule parabolics of superalgebras]{Cominuscule parabolics of simple finite dimensional Lie superalgebras }

\author[Dimitar Grantcharov and Milen Yakimov]{ 
Dimitar Grantcharov$^1$
and
Milen Yakimov$^2$
}
\thanks{$^1$The author 
was supported in part by NSA grant H98230-10-1-0207 and by Max Planck Institute for Mathematics, Bonn.} 
\thanks{$^2$The author was supported in part
by NSF grant DMS-1001632.}

\date{}
\address{Department of Mathematics, University of Texas at Arlington,
Arlington, TX 76019, U.S.A.} 
\email{grandim@uta.edu}
\address{Department of Mathematics, Louisiana State University,
Baton Rouge, LA 70803, U.S.A.}
\email{yakimov@math.lsu.edu} 
\keywords{Simple Lie superalgebras, cominuscule parabolic subalgebras, 
Levi subalgebras, nilradicals}
\subjclass[2010]{Primary 17B05; Secondary 17B22}
\begin{abstract} 
We give an explicit classification 
of the cominuscule parabolic subalgebras of all 
complex simple finite dimensional Lie superalgebras.
\end{abstract}
\maketitle
\label{Intro}
Cominuscule parabolic subalgebras of simple finite dimensional 
Lie algebras play an important role in representation theory, geometry,
and combinatorics.
In the complex case the corresponding 
flag varieties are Hermitian symmetric spaces of compact type
\cite{H}. There has been a lot of research on properties 
of  Schubert varieties in those spaces, see Chapter 9 of \cite{BL} for a comprehensive
survey. More recently it was proved that many results for 
Schubert calculus on Grassmannians extend to those \cite{PSo, TY}
and that the totally nonnegative part of cominuscule 
flag varieties can be described much more explicitly 
than the general case 
in terms of Le diagrams \cite{P,LW}. The standard Poisson 
structure on cominuscule flag varieties also has special properties: 
the (finitely many) orbits of the standard Levi subgroup are complete Poisson
submanifolds \cite{BGY,GY} and each of them is a quotient of the 
standard Poisson structure of the Levi subgroup.
Finally, one should note that cominuscule parabolics
have special properties in numerous other respects.

Parabolic subalgebras of Lie superalgebras are much 
less well understood.
In this paper we address the question of classification of all cominuscule 
parabolic subalgebras of the complex simple finite dimensional Lie superalgebras. 
One of the many equivalent definitions of a cominuscule 
parabolic subalgebra of a simple 
finite dimensional Lie algebra is that it is a parabolic subalgebra with 
abelian nilradical. Already the case of $\gs \gl(m|n)$ presents 
an intrinsic example of parabolic subalgebras, which should be considered 
cominuscule in the super sense. The standard $\Z$-grading of 
$\gs \gl(m|n)$ gives rise to two maximal parabolic subalgebras 
whose (pure odd) ``nilradicals'' are abelian.

In order to define cominuscule parabolics in the super case,
we first reexamine the definition of a parabolic subalgebra 
and nilradical. Unlike the classical (even) case there is no uniform 
definition of either of the two notions in terms of Borel subalgebras
and maximal nilpotent ideals, respectively. For the first one we 
work with root subalgebras. Let $\gg$ be a complex 
simple finite dimensional Lie superalgebra and $\gh$ be a fixed Cartan 
subalgebra. Denote by $\Delta$ the set of roots of $\gg$ with 
respect to $\gh$.
For $\alpha \in \Delta$ let $\gg^\alpha$ be the corresponding 
root space. We call a subalgebra $\gl$ of $\gg$ a 
{\it{root subalgebra}} if it has the form 
\begin{equation}
\label{root}
\gl = (\gl \cap \gh) 
\oplus \left( \bigoplus_{\alpha \in \Phi} \, \gg^\alpha \right)
\end{equation}
for some subset $\Phi \subseteq \Delta$. If $\Delta$ is symmetric
(i.e. $\Delta = - \Delta$), 
then we call a proper subset $P$ of $\Delta$ 
a {\em{parabolic set of roots}} if
\[
\Delta = P  \cup(-P) \; \; \; \mbox{and} \; \; 
\alpha, \beta \in P \; {\text { with }} \; \alpha + \beta \in \Delta 
\; \; 
{\text { implies }} \; \; \alpha + \beta \in P.
\]
If $\Delta \neq - \Delta$, then $P \subsetneq \Delta$ will be called 
parabolic if $P = \widetilde{P} \cap \Delta$ 
for some parabolic subset $\widetilde{P}$ of $\Delta \cup (- \Delta)$. 
We will call a subalgebra of $\gg$ {\em{parabolic}} if it is 
a root subalgebra as in \eqref{root} for a parabolic 
subset of roots $\Phi \subsetneq \Delta$ and contains 
$\gh$. In other words, given a parabolic subset of roots $P$,
then the corresponding parabolic subalgebra of $\gg$ is 
\[
\gp_P := \gh \oplus \left( \bigoplus_{\alpha \in P} \, \gg^\alpha \right).
\]

For a symmetric root systems $\Delta$ and a parabolic set of roots 
$P \subsetneq \Delta$, we call $L := P \cap (-P)$ 
the {\em{Levi component}} of $P$, $N^+ := P \backslash (-P)$ the 
{\em{nilradical}} of $P$, and $P = L \sqcup N^+$ the Levi decomposition 
of $P$. In the nonsymmetric case ($\Delta \neq - \Delta$)
one cannot use the same
formulas, since $N^+ = P \backslash (-P)$ is generally
not closed under addition (i.e. $\alpha, \beta \in N^+$, 
$\alpha + \beta \in \Delta$ does not imply $\alpha + \beta \in N^+$).    
If $\Delta \neq - \Delta$, then we choose 
a parabolic subset $\widetilde{P}$ of $\Delta \cap (-\Delta)$ such that
$P = \widetilde{P} \cap \Delta$, and define
\[
\widetilde{L} = \widetilde{P} \cap ( - \widetilde{P}) \; \;
\mbox{and} \; \;  
\widetilde{N}^+ = \widetilde{P} \backslash (- \widetilde{P} ).
\]
We call $L : = \widetilde{L} \cap P$ a Levi component of $P$, 
$N^+ = \widetilde{N}^+ \cap P$ a nilradical of $P$, 
and $P = L \sqcup N^+$ a Levi decomposition of $P$.
We note that in the nonsymmetric case
the definition of a Levi component and nilradical 
of $P$ essentially depends on the 
choice of a parabolic subset $\widetilde{P}$ of 
$\Delta \cap (-\Delta)$.
We refer the reader to Remarks \ref{nonu1} and \ref{nonu2}
for details.

Let $P$ be a parabolic subset of roots of $\gg$. For a 
Levi decomposition $P = L \sqcup N^+$ of $P$ we define 
the subalgebras
\[
\gl := \gh \oplus \left( \bigoplus_{\alpha \in L} \, \gg^\alpha \right) \; \;
\mbox{and}
\; \; 
\gn^+ := \bigoplus_{\alpha \in N^+} \, \gg^\alpha
\]
of $\gp_P$, and call them a {\em{Levi subalgebra}} 
and {\em{nilradical}} of $\gp_P$. This gives rise to the 
semidirect sum decomposition $\gg_P = \gl \ltimes \gn^+$,
which will be called a {\em{Levi decomposition}} of $\gp_P$.
(Here and below the symbol $\ltimes$ will stand for 
semi-direct sums of Lie superalgebras.) 

We call a parabolic subalgebra $\gp_P$ of $\gg$ 
{\em{cominuscule}} if it has a nilradical $\gn^+$, 
which is abelian. In this paper we investigate 
the parabolic subalgebras of all complex simple 
finite dimensional Lie superalgebras $\gg$. On a case 
by case basis we classify all of their cominuscule 
parabolics. One remarkable consequence of our classification 
is:
 
\begin{theorem} 
\label{unique}
Each cominuscule parabolic subalgebra of a 
complex simple finite dimensional Lie superalgebra 
has a unique Levi decomposition. 
\end{theorem} 

For each cominuscule parabolic subalgebra we describe 
explicitly its Levi subalgebra $\gl$ and the structure of its 
nilradical $\gn^+$ considered as an $\gl$-module. 

We call a parabolic set of roots of $\gg$ cominuscule 
if $\gp_P$ is a cominuscule parabolic subalgebra of $\gg$. 
In the classical even case, one approach to 
comuniscule parabolic subalgebras is through the properties 
of their root systems. Our treatment is based on the 
following super version of this approach, which we prove in 
Proposition \ref{roots}: 

{\em{If $\gg$ is a simple 
finite dimensional Lie superalgebra and 
$\gg \neq S(n), \gg \neq S'(n)$, $\gg \neq \gp \gs \gl (3|3)$,
then a parabolic subset $P$ of $\Delta$ is cominuscule if and only 
if it has a nilradical $N^+$ such that for every $\alpha, \beta$ in $N^+$, 
$\alpha + \beta \notin \Delta$.}} 

This is deduced from the 
fact that for all 
$\gg \neq S(n), \gg \neq S'(n)$, $\gg \neq \gp \gs \gl (3|3)$:
\[
\mbox{if} \; \;  \alpha, \beta, \alpha + \beta \in \Delta,
\; \; \mbox{then} \; \;  
[\gg^\alpha, \gg^\beta] \neq 0.
\]
(In the super case the root spaces $\gg^\alpha$ can have 
dimension more than $1$ and generally
$[\gg^\alpha, \gg^\beta] \neq
\gg^{\alpha + \beta}$.) 
We classify the cominuscule parabolic subsets $P$ of 
$\Delta$ using the above result and 
combinatorial arguments for root systems. 
We prove that a version of the above result is still valid
for $\gg = \gp \gs \gl (3|3)$ and $\gg =S(n), S'(n)$.
In this case
one has to identify the root systems of $\gp \gs \gl (3|3)$
and $S(n), S'(n)$ with subsets of the root systems of $\gs \gl (3|3)$ 
and $W(n)$, respectively, and compute the sums 
of roots $\alpha + \beta$ in the root systems 
of the latter family of Lie superalgebras
(see \S \ref{psl} and \S \ref{S_n} for details).

For each simple finite dimensional Lie superalgebra $\gg$ with 
root system $\Delta$ 
there is a canonical Weyl
group that acts on the set of its parabolic subsets of roots and permutes
the cominuscule parabolic subsets of $\Delta$. We classify the cominuscule 
parabolic subsets of roots of $\gg$ up to the action of this 
Weyl group. For the different types of simple Lie superalgebras
the Weyl group is described as follows. If $\gg$ is a classical 
Lie superalgebra, then $\gg_{\bar{0}}$ is a reductive 
Lie algebra. Its Weyl group $W_{\gg_{\bar{0}}}$ acts in a canonical way on 
the set of the parabolic subsets of $\Delta$ and permutes
the cominuscule parabolic subsets of $\Delta$. 
We classify the cominuscule parabolics of the {\em{basic 
classical Lie superalgebras}} and the {\em{strange classical 
Lie superalgebras}} up to the action of this group
in Sections \ref{classic} and \ref{strange}, 
respectively. The second family of Lie superalgebras 
contains nonsymmetric root systems $\Delta$.
Those are the first ones for which we establish the validity 
of Theorem \ref{unique}. Finally the classification 
of cominuscule parabolics of the {\em{Cartan type}} Lie superalgebras
is carried out in Section \ref{Cartan}. In this case 
$\gg_{\bar{0}}$ has a canonical Levi subalgebra $\gl_{\bar{0}}$. 
Its Weyl group $W_{\gl_{\bar{0}}}$ acts in a natural way 
on the set of parabolic subsets of $\Delta$ and permutes 
the cominuscule parabolic subsets of $\Delta$. The latter are 
classified up to the action of the Weyl group $W_{\gl_{\bar{0}}}$.

To keep the size of the paper down, we will 
not summarize the results of the classification 
theorems of Sections \ref{classic}, \ref{strange}, 
and \ref{Cartan}. The subsections of those sections 
are labeled by the corresponding simple Lie superalgebras, 
so the reader can easily search those results. Another 
interesting
corollary of our classification is that {\em{all cominuscule 
parabolic sets of roots for simple finite dimensional Lie superalgebras
are principal}}, i.e. they come from {\em{triangular decompositions}} 
of the root systems (see \S \ref{sec_principal} for definitions). 

We expect that the class of cominuscule parabolics 
of simple Lie superalgebras will play an important 
role, similar to the ones in the even case.  
In particular, we expect that parabolic induction 
from such will behave well and that the cominuscule 
super flag varieties will have many special properties 
distinguishing them from general super flag varieties.

We finish the introduction with several notational conventions, which will 
be used throughout the paper. We will denote the standard representations 
of $\gg \gl (p)$ and $\gg \gl (p|q)$ by $V^p$ and $V^{p|q}$, 
respectively. The same notation will be used for 
the restrictions of these representations to the subalgebras
of $\gg \gl (p)$ and $\gg \gl (p|q)$. We will use $S^k$ and 
$\bigwedge^k$ to denote the $k$th (super)symmetric power and 
(super)exterior power, respectively. For a module $M$, 
$M^*$ will stand for its dual module. We will follow \cite{P} in our 
notation for Lie superalgebras, except that we will denote by 
$S'(n)$ the Lie superalgebras series denoted by $\widetilde{S}(n)$ 
in \cite{P}. For a Lie superalgebra $\ga$, $\ga'$ will stand for
its derived subalgebra $[\ga,\ga]$. Set-theoretic unions will 
be denoted by $\cup$ and disjoint unions will be denoted 
by $\sqcup$.

\section{Parabolic sets of roots and parabolic subalgebras}
\setcounter{equation}{0}

\label{Prelim}
This section contains some general facts about 
parabolic sets of roots and parabolic subalgebras 
of simple finite dimensional 
Lie superalgebras $\gg$. We define Levi components and nilradicals
of parabolic sets of roots, 
and use those to define Levi subalgebras and nilradicals 
of parabolic subalgebras of $\gg$. In the special 
case of principal parabolic sets of roots those recover the
triangular decompositions of $\gg$. We define 
cominuscule parabolic subalgebras and establish 
a relationship to the properties of the related 
parabolic sets of roots.

\subsection{Levi decompositions of parabolic sets of roots.}
\label{LN}
In what follows, unless otherwise stated $\gg = \gg_{\bar{0}} \oplus \gg_{\bar{1}}$ will
denote a simple finite dimensional Lie superalgebra  over $\C$ (see
\cite{K} and \cite{Sch} for details). A {\it Cartan subalgebra}
$\gh = \gh_{\bar{0}} \oplus \gh_{\bar{1}}$ of $\gg$, is by
definition a selfnormalizing nilpotent subalgebra. Then
$\gh_{\bar{0}}$ is a Cartan subalgebra of
$\gg_{\bar{0}},$ and $\gh_{\bar{1}}$ is the maximal subspace of
$\gg_{\bar{1}}$ on which $\gh_{\bar{0}}$ acts nilpotently (see
\cite[Proposition 1]{PS} for a proof). We denote by $\Delta =
\Delta{(\gg, \gh)}$ the {\it roots of $\gg$ with respect to $\gh$}.
Thus $\Delta = \{ \alpha \in \gh_{\bar{0}}^*, \; \alpha \neq 0 \; |
\; \gg^{\alpha} \neq 0\}$. For $\bar{\imath} \in \Z/2 \Z$ we set
$\Delta_{\bar{\imath}} = \{ \alpha \in \Delta\; | \;
\gg_{\bar{\imath}}^{\alpha} \neq 0\}$.

\begin{definition}
\label{par1} Let $\Delta$ be the root system of a simple 
finite dimensional Lie superalgebra.
If $\Delta = - \Delta$, we will call a proper subset $P$ of $\Delta$ 
a {\em parabolic set of roots} if
\begin{equation} \label{P1}
\Delta = P \cup (-P) \; \; \; \mbox{and} 
\end{equation}
\begin{equation} \label{P2}
\alpha, \beta \in P \; {\text { with }} \; \alpha + \beta \in \Delta 
\; \; 
{\text { implies }} \; \; \alpha + \beta \in P.
\end{equation}

If $\Delta \neq - \Delta$, $P \subsetneq \Delta$ will be called 
a parabolic subset if $P = \widetilde{P} \cap \Delta$ 
for some parabolic subset $\widetilde{P}$ of $\Delta \cup (- \Delta)$. 
\end{definition}

Next we define a Levi decomposition of a parabolic set of roots.

\begin{definition}
\label{dLN} Let $P$ be a parabolic set of roots of the root system 
$\Delta$ of a simple finite dimensional Lie superalgebra. 

(1) If $\Delta = - \Delta$ we will call $L := P \cap (-P)$ 
{\em the Levi component} of $P$, $N^+ := P \backslash (-P)$ {\em the 
nilradical} of $P$, and $P = L \sqcup N^+$ {\em the Levi decomposition} 
of $P$.  

(2) If $\Delta \neq - \Delta$, then we choose 
a parabolic subset $\widetilde{P}$ of $\Delta \cap (-\Delta)$ such that
$P = \widetilde{P} \cap \Delta$ and set
\[
\widetilde{L} = \widetilde{P} \cap ( - \widetilde{P}), \quad 
\widetilde{N}^+ = \widetilde{P} \backslash (- \widetilde{P} ).
\]
We call $L : = \widetilde{L} \cap P$ {\em a Levi component} of $P$, 
$N^+ = \widetilde{N}^+ \cap P$ {\em a nilradical} of $P$, 
and $P = L \sqcup N^+$ {\em a Levi decomposition} of $P$.
\end{definition}

In the nonsymmetric case $\Delta \neq - \Delta$, 
the definition of a Levi component and nilradical 
of $P$ essentially depends on the 
choice of a parabolic subset $\widetilde{P}$ of 
$\Delta \cap (-\Delta)$ such that $P = \widetilde{P} \cap \Delta$.
We provide examples and discuss this further 
in Remarks \ref{nonu1} and \ref{nonu2}.

The following lemma contains several simple facts 
for the Levi components and nilradicals of parabolic sets of roots.

\begin{lemma} \label{sums}
Let $P$ be a parabolic subset of $\Delta$ and let $P = \widetilde{P} \cap \Delta$ 
for  some parabolic subset $\widetilde{P}$ of $\Delta \cup (- \Delta)$. 
Set $\widetilde{L} := \widetilde{P}  \cap (- \widetilde{P})$ and
$\widetilde{N}^{\pm} := ( \pm \widetilde{P})  \setminus (\mp \widetilde{P} )$. Let 
$L:=\widetilde{L} \cap \Delta $ and $N^{\pm}:=\widetilde{N}^{\pm} \cap \Delta $
be the Levi component and nilradical of $P$ corresponding to $\widetilde{P}$.
Then 
\[
\Delta \cup (- \Delta) = \widetilde{N}^{-} \sqcup \widetilde{L} 
\sqcup \widetilde{N}^{+}, \quad
\Delta = N^{-} \sqcup L \sqcup N^{+},
\]
and
\[
\widetilde{L} = L \cup (-L), \quad
\widetilde{N}^{\pm } = N^{\pm} \cup (- N^{\mp}).
\] 
If $\alpha_L, \alpha_L' \in L$ and 
$\alpha_{N^{\pm}}, \alpha'_{N^{\pm}} \in N^{\pm}$ then

(i) $- \alpha_{N^{\pm}}  \in \Delta$ implies 
$- \alpha_{N^\pm} \in N^{\mp}$;

(ii) $\alpha_L + \alpha_{N^{\pm}} \in \Delta$ 
implies $\alpha_L + \alpha_{N^\pm} \in N^{\pm}$;

(iii) $\alpha_L + \alpha_L' \in \Delta $ implies $\alpha_L + \alpha_L' \in L$;

(iv) $\alpha_{N^{\pm}} + \alpha_{N^{\pm}}' \in \Delta$ implies $\alpha_{N^{\pm}} + \alpha_{N^{\pm}}' \in N^{\pm}$.

\end{lemma}
\noindent
{\bf Proof.} The set theoretic identities are easy to deduce and are left to 
the reader.

Since $P^- = L \sqcup N^{-} = (\widetilde{L} \sqcup \widetilde{N}^{-}) \cap \Delta$
is a parabolic subset of $\Delta$, it is sufficient to prove (i)--(iv) 
for $L$ and $N^+$. From the definition of $L$ and $N^+$ it also 
follows that it is enough to consider the case of 
$\Delta = - \Delta$, i.e. $\widetilde{P} = P$. 
 
For (i), if $- \alpha_{N^{+}} \in P$ then $\alpha_{N^{+}} \in P \cap (-P)$, 
which is a contradiction.  Assume that $\alpha_L + \alpha_{N^{+}} \in L$ in (ii). 
Hence $- \alpha_L - \alpha_{N^{+}} \in L$ which together with 
$\alpha_L \in L$ imply $\pm \alpha_{N^{+}} \in P$ 
leading to a contradiction. 
If $\alpha_L + \alpha'_L \in N^{+}$ in (iii), then by (ii), 
$\alpha_L = - \alpha'_L + (\alpha_L + \alpha'_L) \in N^{+}$, 
which is a  contradiction.  Finally, for (iv), assume that 
$\alpha_{N^{+}} + \alpha_{N^{+}}' \in L$. 
Then $- \alpha_{N^{+}} - \alpha_{N^{+}}' \in L$ and $ \alpha_{N^{+}} \in N^{+}$ 
imply $- \alpha_{N^{+}}' \in L$ by (ii), and again we 
reach a contradiction. \hfill $\square$
\subsection{Principal parabolic sets of roots.} \label{sec_principal}
\label{principal}
Let $V$ be a finite dimensional real vector space such that 
$\Delta \subset V \setminus \{0\}$. A partition 
$\Delta = \Delta^- \sqcup \Delta^0 \sqcup \Delta^+$ 
is called  a {\it triangular decomposition} of $\Delta$ 
if there exists a functional $\Lambda \in V^*$ such that 
$\Delta^0 = \Delta \cap \mbox{Ker} \Lambda$ and
$\Delta^\pm = \{ \alpha \in \Delta \, | \, \Lambda(\alpha) \gtrless 0 \}$. 
A subset $P$ of $\Delta$ is called a {\it principal parabolic set} 
if there exists a triangular decomposition
$\Delta = \Delta^- \sqcup \Delta^0 \sqcup \Delta^+$ such that
$P = \Delta^0 \sqcup \Delta^+$. In such a case we
write $P = P(\Lambda)$. 

\begin{proposition} 
\label{prop1}
Every principal parabolic subset $P$ of the set of roots $\Delta$ 
of a simple finite dimensional Lie superalgebra
is a parabolic subset of $\Delta$. The set 
$\Delta^0$ is a Levi component of $P$ and 
$\Delta^+$ is a nilradical of $P$.
\end{proposition}
\noindent
{\bf Proof.}
Consider the principal parabolic subset of $\Delta \cup (-\Delta)$:
\[
\widetilde{P} = \{ \alpha \in \Delta \cup (- \Delta)  \, | \, 
\Lambda(\alpha) \geq 0 \}.
\]
It is clear that $\widetilde{P}$ is a parabolic subset of $\Delta \cup (-\Delta)$
and that $P = \widetilde{P} \cap \Delta$. This implies the 
first statement. For the second statement, observe that
\[
\widetilde{L} = \widetilde{P} \cap (- \widetilde{P}) = 
\{ \alpha \in \Delta \cup (- \Delta)  \, | \, 
\Lambda(\alpha) = 0 \} = \Delta^0 \cup (- \Delta^0)
\]
and 
\[
\widetilde{N}^+ = \widetilde{P} \backslash (- \widetilde{P}) = 
 \{ \alpha \in \Delta \cup (- \Delta)  \, | \, 
\Lambda(\alpha) > 0 \} =
\Delta^+ \cup (- \Delta^-).
\]
Therefore $L = \widetilde{L} \cap \Delta = 
\{ \alpha \in \Delta \cup (- \Delta)  \, | \, 
\Lambda(\alpha) = 0 \} \cap \Delta = \Delta^0$
and $N^+ = \widetilde{N}^+ \cap \Delta = \{ \alpha \in \Delta \cup (- \Delta)  \, | \, 
\Lambda(\alpha) > 0 \} \cap \Delta = \Delta^+$.
\hfill $\square$

\begin{remark}
\label{nonu1} A principal parabolic subset $P$ of $\Delta$ can have 
different Levi decompositions when $\Delta \neq - \Delta$. Given 
$P$, define the polyhedron
\begin{equation}
\label{FP}
{\mathcal{F}}(P) = \{ 
\Lambda \in V^* \; | \; \Lambda(\alpha) \geq 0, 
\, \forall \alpha \in P, \; \Lambda(\alpha) < 0,
\, \forall \alpha \in \Delta \backslash P \}.
\end{equation}
(Here and below we use the term polyhedron in the wide sense 
as a subset of a real vector space, which is an intersection 
of a finite collection of open and closed half spaces.) 
There are some immediate consequences of the inequalities 
in \eqref{FP}. For instance, the first condition in \eqref{FP} 
implies $\Lambda (\alpha) = 0$, $\forall \alpha \in P \cap (-P)$.
For each functional $\Lambda \in {\mathcal{F}}(P)$ define 
\[
L(\Lambda) = \{ \alpha \in \Delta \; | \; 
\Lambda (\alpha) = 0 \}, \; \; 
N^+(\Lambda) = \{ \alpha \in \Delta \; | \; 
\Lambda (\alpha) > 0 \}.
\]
It follows from Proposition \ref{prop1} 
that $P = L(\Lambda) \sqcup N^+(\Lambda)$ 
is a Levi decomposition of $P$, 
$\forall \Lambda \in {\mathcal{F}}(P)$. 
This Levi decomposition is the same for points $\Lambda$ in the interior 
of a fixed face of ${\mathcal{F}}(P)$, but differs 
for points $\Lambda$ that belong to the interiors 
of different faces of ${\mathcal{F}}(P)$. 
This is further illustrated in Remark \ref{nonu2}. 
\end{remark}

The converse to the first statement of Proposition \ref{prop1}
is true for finite dimensional reductive Lie algebras
(see, for example,  \cite[Proposition VI.7.20]{Bo}).
More generally, we have (see \cite[Proposition 2.10]{DFG}):

\begin{proposition} 
\label{kac-moody} 
Let $\gg$ be a quasisimple regular Kac-Moody superalgebras 
and let $P$ be a parabolic subset of $\Delta$. Then  $P$ is a 
principal parabolic subset of $\Delta$.
\end{proposition}

The simple finite dimensional Lie superalgebras which are not Lie algebras are: 
$\gs \gl (m|n)$ for $m \neq n$, $\gp \gs \gl(m|m)$, $\go \gs \gp(m|2n)$,
$D(2,1;\alpha)$, $F(4)$, $G(3)$, ${\bf s} \gp (n)$, $\gp \gs \gq(n)$, and the 
Cartan type superalgebras $W(n)$, $S(n)$, $S'(n)$, and $H(n)$. For the
restrictions on the parameters $m,n$, and $\alpha$ as well as isomorphisms 
among the superalgebras listed above we refer the reader to \cite{K}.
Among those, the quasisimple regular Kac--Moody Lie superalgebras are 
$\gs \gl(m|n)$ for $m \neq n$,
$\go \gs \gp(m|2n)$, $D(2,1;\alpha)$, $F(4)$, and $G(3)$, 
(see \cite{S2}). Proposition \ref{kac-moody} applies to 
them, i.e. all parabolic sets of roots for them are principal.

\subsection{Cominuscule parabolic subalgebras.}
\label{cominuscule}
Recall from \S \ref{LN} that $\gh$ denotes a fixed Cartan subalgebra
of the simple finite dimensional Lie superalgebra $\gg$ and that 
$\Delta = \Delta(\gg, \gh)$ denotes the corresponding set of roots of $\gg$. 
For $\alpha \in \Delta$, we will denote by $\gg^\alpha$ the corresponding 
root space. We will call a subalgebra of $\gg$ a 
{\it{root subalgebra}} if it has the form 
\begin{equation}
\label{root-sub}
\gl = (\gl \cap \gh) 
\oplus \left( \bigoplus_{\alpha \in \Phi} \, \gg^\alpha \right)
\end{equation}
for some subset $\Phi \subseteq \Delta$. We will call 
a root subalgebra of $\gg$ {\it{parabolic}} if 
$\gl \supseteq \gh$ and 
$\Phi$ is a parabolic subset of $\Delta$.
\begin{definition}
\label{LN-alg}
Let $P$ be a parabolic subset of the set of roots 
$\Delta$ of a simple finite dimensional Lie superalgebra $\gg$
with respect to a Cartan subalgebra $\gh$ and  
\begin{equation} 
\label{parabolic}
\gp_P = \gh \oplus \left( \bigoplus_{\alpha \in P} \, \gg^\alpha \right)
\end{equation} 
be the corresponding parabolic subalgebra. Given a  
Levi decomposition $P = L \sqcup N^+$ of $P$ we define 
the subalgebras
\begin{equation}
\label{nn}
\gl = \gh \oplus \left( \bigoplus_{\alpha \in L} \, \gg^\alpha \right) \; \;
\mbox{and}
\; \; 
\gn^+ = \bigoplus_{\alpha \in N^+} \, \gg^\alpha
\end{equation}
of $\gp_P$, and call them {\em a Levi subalgebra} and {\em a nilradical} of $\gp_P$, respectively. 
The semidirect sum decomposition $\gp_P = \gl \ltimes \gn^+$
will be called {\em a Levi decomposition} of $\gp_P$.
\end{definition}

In the classical even case ($\gg_{\bar{1}}=0$) 
the root system is symmetric $\Delta = - \Delta$ 
and the above decomposition $\gp_P = \gl \ltimes \gn^+$
is precisely the Levi decomposition of $\gp_P$
for the unique Levi subalgebra $\gl$ containing 
$\gh$.

The next definition singles out the class of cominuscule parabolic
subalgebras.
\begin{definition}
\label{comin}
We call a parabolic subalgebra of a complex simple 
finite dimensional Lie superalgebra $\gg$ {\em cominuscule},
if has a nilradical which is abelian. A parabolic subset $P$
of the set of roots $\Delta$ of $\gg$ will be called 
cominuscule if the corresponding 
parabolic subalgebra $\gg_P$ is cominuscule. 
\end{definition}

If $\gg$ is a complex simple finite dimensional Lie algebra,
then Definition \ref{comin} singles out exactly the 
class of cominuscule parabolic subalgebras of $\gg$ which 
contain the fixed Cartan subalgebra $\gh$.

In the Appendix we will prove the following proposition.

\begin{proposition} 
\label{appendix2}
Let $\gg$ be a simple 
finite dimensional Lie superalgebra, 
$\gg \neq S(n), \gg \neq S'(n)$, $\gg \neq \gp \gs \gl (3|3)$, 
$\gh$ be a Cartan subalgebra of $\gg$, 
and $\Delta$ be the corresponding root system.
If $\alpha, \beta, \alpha + \beta \in \Delta$, 
then 
\[
[\gg^\alpha, \gg^\beta] \neq 0.
\]
\end{proposition}
We note that in the super case, 
$\alpha, \beta, \alpha + \beta \in \Delta$
does not imply
\[
[\gg^\alpha, \gg^\beta] =
\gg^{\alpha + \beta}.
\]

The next proposition generalizes an important 
property of cominuscule parabolic subalgebras
which holds in the classical (even) case. It reduces 
the problem of classification of cominuscule parabolics
of simple finite dimensional Lie superalgebras to a problem 
for the corresponding root systems and will be extensively used in 
the paper.

\begin{proposition} 
\label{roots}
Let $\gg$ be as in Proposition \ref{appendix2}. A parabolic subset $P$ 
of the set of roots of $\gg$ is cominuscule if and only if it has 
a nilradical $N^+$ such that for every $\alpha, \beta$ in $N^+$, 
$\alpha + \beta \notin \Delta$,
\end{proposition}
\noindent
{\bf Proof.} If a parabolic subset of roots has the above property,
then for all $\alpha, \beta \in N^+$, $[\gg^\alpha, \gg^\beta]= 0$.
Therefore the nilradical $\gn^+$ of the parabolic subalgebra $\gp_P$
corresponding to $N^+$ is abelian, cf. \eqref{parabolic} and \eqref{nn}, 
and thus $\gp_P$ is cominuscule.

In the other direction, assume that the parabolic subset 
of roots $P$ is cominuscule. Let $N^+$ be a nilradical of $P$ 
such that the corresponding nilradical $\gn^+$ of $\gg$ given by 
\eqref{nn} is abelian. 
If there exist $\alpha, \beta \in N^+$
such that $\alpha + \beta \in \Delta$, then by 
Proposition \ref{appendix2}, $[\gg^\alpha, \gg^\beta] \neq 0$. 
This is a contradiction since $\gg^\alpha$ and $\gg^\beta$ 
are subspaces of the nilradical of the parabolic 
subalgebra $\gg_P$. 
Therefore for all $\alpha, \beta$ in $N^+$, 
$\alpha + \beta \notin \Delta$, which completes the 
proof of the proposition.
\hfill $\square$

\begin{remark}
The cases $\gg = \gp \gs \gl (3|3)$, $\gg = S(n)$ and $\gg = S'(n)$ 
require special attention and modifications of Proposition \ref{appendix2} 
in these cases are established in \S \ref{psl} and \S \ref{S_n}, respectively.
\end{remark}

The classification of all cominuscule parabolic subsets established 
in the next sections implies the following result. 
\begin{theorem}
All cominuscule parabolic subsets of the simple finite dimensional 
Lie superalgebras are principal.
\end{theorem}
\subsection{Passing to subalgebras.}
\label{even}
We will need the following lemma for reduction 
of cominuscule parabolics to certain subalgebras. 
Its proof is straightforward and will be left to the reader.
\begin{lemma} \label{restrict} Let $\gg$ be a simple Lie superalgebra 
and $\gh$ be a Cartan subalgebra of $\gg$. Let
$\ga$ be subalgebra of $\gg$, which is either 
a simple superalgebra or an (even) reductive Lie algebra.
Assume that $\ga \cap \gh$ is a Cartan subalgebra of $\ga$
and that $\alpha|_{\ga \cap \gh_{\bar{0}}} \neq 
\beta|_{\ga \cap \gh_{\bar{0}}}$ for all roots $\alpha \neq \beta$ 
of $\gg$.
Let $\Delta$ be the root system of $\gg$ with respect to $\gh$
and $\Delta_{\ga}$ be the root system of $\ga$ 
(considered as a subset of $\Delta$) with respect to 
$\ga \cap \gh$. 

If $P$ is a parabolic subset of $\Delta$, then 
$P \cap \Delta_{\ga}$ is either equal to $\Delta_{\ga}$ 
or to a parabolic subset 
of $\Delta_{\ga}$. In the latter case, if $P = L \sqcup N^+$ is a 
Levi decomposition of $P$, then 
$P \cap \Delta_{\ga}= (L \cap \Delta_{\ga}) \sqcup (N^+ \cap \Delta_{\ga})$ 
is a Levi decomposition of $P \cap \Delta_{\ga}$.
\end{lemma}

Finally, we will make extensive use of the classification of cominuscule 
parabolics of classical simple Lie algebras. We recall it below 
for completeness. We will use the notation of \cite{Bo}
for the root spaces of $\gs \gl (n)$, $\gs \go (2n)$, 
$\gs \go (2n+1)$, and $\gs \gp (2n)$.

\begin{proposition} \label{algebras}
The following list describes the Levi components and nilradicals of all cominuscule 
parabolic sets of roots for the finite dimensional simple Lie algebras $\gg$ 
of type $A,B,C,$ and $D$ up to the action of the Weyl group $W_{\gg}$ of $\gg$.

(i) $L_{\gs \gl (n)} (n_0)= \{ \varepsilon_i - \varepsilon_j  \; | 
\; i \neq j \leq n_0 \mbox{ or } i \neq j > n_0 \}$ and 
$N_{\gs \gl (n)}^+ (n_0)= \{ \varepsilon_i - \varepsilon_j  
\; | \; i \leq n_0 < j \},$ $1\leq n_0 \leq n-1$, 
for $\gg = \gs \gl (n)$.  

(ii) $L_{\gs \go (2n+1)} = \{ \pm \varepsilon_i \pm  \varepsilon_j,  
\pm \varepsilon_i   \; | \; j>i>1 \}$ and $N_{\gs \go (2n+1)}^+ = 
\{ \varepsilon_1 \pm  \varepsilon_j, 
\varepsilon_1 \; | \; j > 1\}$ for $\gg = \gs \go (2n+1)$.

(iii) $L_{\gs \gp (2n)} = \{  \varepsilon_k -  \varepsilon_l  \; | \; k \neq l \}$ 
and $N_{\gs \gp (2n)}^+ = \{   \varepsilon_k +  \varepsilon_l, 2  \varepsilon_k  \; | 
\; k \neq l \}$ for $\gg = \gs \gp (2n)$.

(iv)  For $\gg = \gs \go (2n)$, $n \geq 2$:
\begin{eqnarray*}  
&&L_{\gs \go (2n)} (1)=   \{ \pm \varepsilon_i 
\pm  \varepsilon_j,  \; | \; j>i>1 \} \; \; \mbox{ and } 
\; \; N_{\gs \go (2n)}^+(1) =  
\{ \varepsilon_1 \pm  \varepsilon_j \; | \; j > 1\},\\
&&L_{\gs \go (2n)} (n) =   \{  \varepsilon_i  -  \varepsilon_j,  \; | \;  i \neq j \} \; \; 
\mbox{ and } \; \; N_{\gs \go (2n)}^+ (n)=  \{ \varepsilon_i +  \varepsilon_j \; | 
\; i \neq j \},\\
&&\overline{L}_{\gs \go (2n)} (n) = (\theta L_{\gs \go (2n)} (n))\; \; 
\mbox{ and } \; \; \overline{N}_{\gs \go (2n)}^+ (n) =  
\theta (N_{\gs \go (2n)}^+ (n)),
\end{eqnarray*}
where $\theta$ is the involutive automorphism of the Dynkin diagram $D_n$ 
that preserves $ \varepsilon_i  -  \varepsilon_{i+1}$, $i <n-1$, and 
interchanges $ \varepsilon_{n-1} -  \varepsilon_{n}$ and 
$ \varepsilon_{n-1} +  \varepsilon_{n}$.
\end{proposition}

For convenience we will also use the notation from 
Proposition \ref{algebras} (i) for $n_0 = n$, i.e. for $\gg = \gs \gl (n)$ 
we set $L_{\gs \gl (n)} (n):=\Delta$ and 
$N_{\gs \gl (n)}^+ (n)=\emptyset$.

\section{Classification of the cominuscule parabolics of basic classical Lie superalgebras}
\label{classic}
\setcounter{equation}{0}

In this section we classify the cominuscule parabolic subsets of all basic classical Lie superalgebras.
Fix such a Lie superalgebra $\gg$ and a Cartan subalgebra $\gh$ of it. Denote, as before, 
the set of roots of $\gg$ with respect to $\gh$ by $\Delta$. For this 
class of superalgebras, the even part $\gg_{\bar{0}}$ of $\gg$ is a reductive 
Lie algebra. The Weyl group $W_{\gg_{\bar{0}}}$ of the even part $\gg_{\bar{0}}$ 
acts on the root system
$\Delta$ of $\gg$ and thus on the parabolic subsets of $\Delta$. 
(Each element of the Weyl group of the even part $w \in W_{\gg_{\bar{0}}}$ can be lifted 
to an automorphism $\sigma_w$ of $\gg$ stabilizing $\gh$, which induces an automorphism 
of the set of roots. The later does not depend on the 
choice of $\sigma_w$ and by abuse of notation will be denoted by $w$.) 
It is obvious that, if a parabolic subset $P$ of $\Delta$ is cominuscule, 
then $w(P)$ is cominuscule for all $w \in W$.
We will classify the cominuscule parabolic subsets of $\Delta$ up to the 
action of the Weyl group $W_{\gg_{\bar{0}}}$.

In all proofs $P$ will denote a cominuscule parabolic set of roots 
in $\Delta$. Since $\Delta$ is symmetric, in the definition of a parabolic subset $P$ 
of $\Delta$ there is no need to consider a parabolic subset $\wt{P}$ (i.e. $\wt{P}=P$), 
and $P$ has a unique Levi decomposition. The latter is given by $P= L \sqcup N^+$,
where the Levi component is $L= P \cap (-P)$ and the nilradical is 
$N^+ = P \backslash (-P)$. Set $P_{\bar{0}}:=P \cap \Delta_{\bar{0}}$, 
$L_{\bar{0}}:=L \cap \Delta_{\bar{0}}$, 
and $N_{\bar{0}}^+:=N^+ \cap \Delta_{\bar{0}}$. By Lemma \ref{restrict}, 
either $P_{\bar{0}} = \Delta_{\bar{0}}$ or  $P_{\bar{0}}$
is a parabolic subset of $\Delta_{\bar{0}}$ with Levi component $L_{\bar{0}}$ and nilradical
$N_{\bar{0}}^+$. As in Definition \ref{LN-alg} we will denote by $\gl$ and $\gn^+$ 
the Levi subalgebra and nilradical of $\gp_P$.
\subsection{$\gg = \gs \gl (m|n)$, $m\neq n$.}
\label{sl} 
In this case $\gg_{\bar{0}} \cong \gs \gl (m) \oplus \gs \gl (n)  \oplus \C$ with
$\Delta_{\bar{0}} = \{ \varepsilon_i - \varepsilon_j, \delta_k - 
\delta_l \; | \; 1 \leq i \neq j \leq m, 1 \leq k \neq l \leq n\}$ and $\Delta_{\bar{1}} 
= \Delta_{\bar{1}}^{+} \sqcup \Delta_{\bar{1}}^{-}$, where $\Delta_{\bar{1}}^{\pm} := 
\{ \pm ( \varepsilon_i - \delta_l) \; | \; 1 \leq i \leq m, 1 \leq l \leq n\}$. 
We introduce the following sets of roots
\begin{eqnarray*}
L_{\gs \gl (m|n)} (m_0 |n_0)& := &  \{ \varepsilon_i - \varepsilon_j \; | \; i, j \leq m_0 
\mbox{ or } i,j > m_0 \} 
\\
& & \sqcup \{ \delta_k - \delta_l,    \; | \; k, l \leq n_0 \mbox{ or } 
k,l > n_0 \} 
\\
& & \sqcup  \{ \pm ( \varepsilon_i - \delta_l)  \; | \; i \leq m_0, l \leq n_0 
\mbox{ or } i >  m_0, l>n_0 \} \mbox{ and }  \\ N^+_{\gs \gl (m|n)} (m_0|n_0) & := & \{ \varepsilon_i - 
\varepsilon_j \; | \; i \leq m_0 <j \} \sqcup \{ \delta_k - \delta_l,    \; | 
\; k \leq n_0<l \} 
\\
& & \sqcup \{\varepsilon_i - \delta_l, \delta_k - \varepsilon_j  \; | 
\; i \leq m_0< j; k \leq n_0 < l\}.
\end{eqnarray*}
Set $P_{\gs \gl (m|n)} (m_0 |n_0):= L_{\gs \gl (m|n)} (m_0 |n_0) \sqcup 
N^+_{\gs \gl (m|n)} (m_0 |n_0)$. 
Note that, if  $m_0=n_0=0$ or $m_0=m,n_0=n$, 
then $P_{\gs \gl (m|n)} (m_0 |n_0)= \Delta$. In all other cases 
$P_{\gs \gl (m|n)} (m_0 |n_0)$ is proper. It is clear that 
$\Delta =  ( - N^+_{\gs \gl (m|n)} (m_0 |n_0) ) \sqcup L_{\gs \gl (m|n)} (m_0 |n_0) \sqcup  
N^+_{\gs \gl (m|n)} (m_0 |n_0)$.

\begin{theorem} \label{prop_slmn}
Every cominuscule parabolic set $P$ of roots of $\gs \gl (m|n)$, $m\neq n$, is conjugated under the action of 
the Weyl group $W_{\gs \gl (m)} \times W_{\gs \gl (n)}$ of $\gg_{\bar{0}}$ to a unique
subset of the form $P_{\gs \gl (m|n)} (m_0 |n_0)$ for some 
$m_0, n_0$, such that $0 \leq m_0 \leq m$,  $0 \leq n_0 \leq n$, and $(m_0, n_0)\neq (0,0), (m,n)$.  
For the Levi subalgebras and nilradicals of the corresponding parabolic subalgebras, 
we have that $\gl \cong \gs \gl (m_0 | n_0 ) \oplus \gs \gl (m-m_0| n-n_0) \oplus \C$ and 
as an $\gl$-module
$\gn^+ \cong V^{m_0|n_0} \otimes \left( V^{m-m_0|n-n_0}\right)^*$. 
\end{theorem}

Recall from the introduction that $V^{p|q}$ denotes the standard 
representation of $\gs \gl(p|q)$.
\\ \hfill \\
\noindent {\bf Proof of Theorem \ref{prop_slmn}.}  Let 
$N^+_{\gs\gl(m)}:=N^+ \cap \{ \varepsilon_i - \varepsilon_j  \; 
| \; 1 \leq i \neq j \leq m \}$ 
and $N^+_{\gs\gl(n)}:=N^+ \cap \{ \delta_k - \delta_l  \; 
| \; 1 \leq k \neq l \leq n \}$. Obviously, 
$N^+_{\bar{0}} = N^+_{\gs\gl(m)} \sqcup N^+_{\gs\gl(n)} $.

By Proposition \ref{kac-moody}, we know that $P= P(\Lambda)$ is a 
principal parabolic set of roots 
determined by some functional $\Lambda$. Let 
$\Lambda (\varepsilon_i - \varepsilon_j) = x_i - x_j,  
\Lambda (\delta_k - \delta_l) = y_k - y_l$, 
and $\Lambda (\varepsilon_i - \delta_k) = x_i - y_k$ 
for some $x_i$ and $y_k$ such that 
$x_1+ \ldots +x_m +y_1+ \ldots +y_n = 0$.

\medskip
\noindent {\it Case 1: 
$N^+_{\gs\gl(m)} \neq \emptyset$ and $N^+_{\gs\gl(n)} \neq \emptyset$}. 
Proposition 
\ref{algebras} implies that $P$ is conjugated under the action of 
the Weyl group $W_{\gs \gl (m)} \times W_{\gs \gl (n)}$
to a unique parabolic subset with
\begin{eqnarray*}
L_{\bar{0}} & = &  \{ \varepsilon_i - \varepsilon_j\; | \; i,j \leq m_0 
\mbox { or } i,j>m_0 \} 
\sqcup \{ \delta_k - \delta_l \; | \;  k,l \leq n_0 \mbox { or } k,l>n_0\}  
\; \; \mbox{and} \\
N^+_{\bar{0}} & = &  \{ \varepsilon_i - \varepsilon_j, \delta_k - \delta_l \; | \; 
i \leq m_0 < j; \; k \leq n_0 <l\}
\end{eqnarray*}
for some $1 \leq m_0 \leq m-1$, $1 \leq n_0 \leq n-1$.
Let us rename $P$ so it has the above even Levi component and nilradical.
Since $\Lambda|_{L_{\bar{0}}} = 0$ and  $\Lambda|_{N^+_{\bar{0}}} $ takes 
positive values we have that 
$x_1= \ldots =x_{m_0} = x$, $x_{m_0+1} = \ldots = x_m = x'$, 
$y_1= \ldots =y_{n_0} = y$, $y_{n_0+1} = \ldots = y_n = y'$ 
for some $x, y, x', y'$ such that $x>x'$, $y>y'$, 
and $m_0 x + (m-m_0)x' + n_0y + (n-n_0)y' = 0$. By Lemma \ref{roots} we have that 
$\varepsilon_{m_0}  - \delta_{n_0}$ is not in $N^+$ and thus $x \leq y$.
Indeed, since 
$\delta_{n_0} - \delta_{n_0 +1} \in N^+$, if   
$\varepsilon_{m_0}  - \delta_{n_0} \in N^+$, it would force 
$\varepsilon_{m_0}-\delta_{n_0+1} = 
(\varepsilon_{m_0}  - \delta_{n_0}) 
+ (\delta_{n_0} - \delta_{n_0 +1}) \in N^+$ which contradicts to $P$ being cominuscule.
On the other hand, since $\delta_{n_0} - \varepsilon_{m_0} \notin N^+$ we have 
$y \leq x$, and consequently $x=y$. Similarly,
$\varepsilon_{m_0 + 1}  - \delta_{n_0+1}$ and $\delta_{n_0+1} - \varepsilon_{m_0+1}$ are not 
in $N^+$ and hence $x' = y'$. The conditions $x = y$, $x' = y'$, and $x > x'$ determine completely 
$P$ and imply that $L = L_{\gs \gl (m|n)} (m_0 | n_0)$ and  
$N^+ = N^+_{\gs \gl (m|n)} (m_0 |n_0)$.

\medskip
\noindent {\it Case 2: $N^+_{\gs\gl(m)} \neq \emptyset$ and 
$N^+_{\gs\gl(n)} = \emptyset$}. 
Like in the previous case, $P$ is conjugated under the 
action of $W_{\gs \gl (m)} \times W_{\gs \gl (n)}$ to a unique parabolic subset
such that
\begin{eqnarray*}
L_{\bar{0}} & = &  \{ \varepsilon_i - \varepsilon_j\; | \; i,j \leq m_0 \mbox { or } 
i,j>m_0 \} \sqcup \{ \delta_k - \delta_l \; | \; 1 \leq  k \neq l \leq n\}  
\; \; \mbox{and} \\
N^+_{\bar{0}} & = &  \{ \varepsilon_i - \varepsilon_j \; | \; i \leq m_0 < j \}
\end{eqnarray*}
for some $1\leq m_0 \leq m-1$. We rename $P$ accordingly.
Then from  $\Lambda|_{L_{\bar{0}}} = 0$ and  $\Lambda|_{N^+_{\bar{0}}} >0$ 
we find $x_1= \ldots =x_{m_0} = x$, $x_{m_0+1} = \ldots = x_m = x'$, 
$y_1= \ldots =y_{n} = y$, 
for some $x, x', y$ such that $x>x'$ and $m_0 x + (m-m_0)x' + ny = 0$. By Lemma 
\ref{roots} we have that $\varepsilon_{m_0+1}  - \delta_{1}$ and 
$\delta_{1} - \varepsilon_{m_0}$ are not in $N^+$. This leads to 
$x \geq y \geq x'$. We proceed with three separate subcases. 
First, if $x=y>x'$ we have
\begin{eqnarray*}
L & = & L_{\bar{0}} \sqcup  \{ \pm ( \varepsilon_i - \delta_l ) \; | \; i \leq m_0, 
1 \leq  l \leq n\}  \; \; \mbox{and} \\
N^+& = &  N^+_{\bar{0}}  \sqcup \{ \delta_k - \varepsilon_j  \; | \;  j>m_0, 
1 \leq  k \leq n \}.
\end{eqnarray*}
In this case $L = L_{\gs \gl (m|n)} (m_0| n)$ and  $N^+ = N^+_{\gs \gl (m|n)} (m_0|n)$, 
i.e. $n_0 = n$.

Second, if $x>y=x'$, then
\begin{eqnarray*}
L & = & L_{\bar{0}} \sqcup  \{ \pm ( \varepsilon_j - \delta_k ) \; | \; j> m_0, l \leq k \leq n\}  \; \; \mbox{and} \\
N^+& = &  N^+_{\bar{0}}  \sqcup \{ \varepsilon_i - \delta_l \; | \;  i \leq m_0, 1 \leq l \leq n\}.
\end{eqnarray*}
This leads to $L = L_{\gs \gl (m|n)} (m_0| 0)$ and  
$N^+ = N^+_{\gs \gl (m|n)} (m_0|0)$, i.e. $n_0 =0$.

And third, if $x>y=x'$, we have that $\varepsilon_1 - \delta_k, 
\delta_k - \varepsilon_{m} \in N^+$ for every $k$, so $P$ is not cominuscule.

\medskip
\noindent {\it Case 3:  $N^+_{\gs\gl(m)} = \emptyset$ and 
$N^+_{\gs\gl(n)} \neq \emptyset$}. Similarly to the previous 
case we verify that $P$ is conjugated under the action of $W_{\gs \gl (m)} \times W_{\gs \gl (n)}$ 
to a unique parabolic subset 
of the form $P_{\gs \gl (m|n)} (0| n_0)$ or $P_{\gs \gl (m|n)} (m| n_0)$.

\medskip
\noindent {\it Case 4: $N^+_{\gs\gl(m)} = \emptyset$ and 
$N^+_{\gs\gl(n)} = \emptyset$}. In this case
\begin{eqnarray*}
L_{\bar{0}} =   \{ \varepsilon_i - \varepsilon_j\; | \; 1 \leq i \neq j \leq m\} \sqcup 
\{ \delta_k - \delta_l \; | \; 1 \leq  k \neq l \leq n\} \; \; \mbox{and} \; \;   
N^+_{\bar{0}} = \emptyset.
\end{eqnarray*}
We have that $x_1= \ldots =x_{m} = x$ and $y_1= \ldots =y_{n} = y$. 
We easily see that $x \neq y$, 
which leads to $L = L_{\bar{0}}$ and  $\gl = \gg_{\bar{0}}$. 
Depending on whether $x>y$ or $x<y$, we obtain 
$N^+ = \{ \varepsilon_i - \delta_l \; | \;  1 \leq i \leq m, 1 \leq l \leq n\}$ 
or $N^+ = \{ \delta_k - \varepsilon_j \; | \;  1 \leq j \leq m, 1 \leq k \leq n\}$. 
These cases correspond to $m_0 = m, n_0 = 0$, and $m_0 = 0, n_0 = n$, respectively.
\hfill $\square$

\begin{remark}
The cases $m_0 = 0, n_0=n$ and $m_0 = m, n_0 = 0$ correspond
to the ``standard'' triangular decomposition 
$\Delta = \Delta_{\bar{1}}^{-}  \sqcup \Delta_{\bar{0}}  
\sqcup\Delta_{\bar{1}}^{+}$, namely  to $P = \Delta_{\bar{0}}  
\sqcup\Delta_{\bar{1}}^{+}$ and $P = \Delta_{\bar{0}}  
\sqcup\Delta_{\bar{1}}^{-}$, respectively.
\end{remark}
\subsection{$\gg = \gp \gs \gl (n|n)$.}
\label{psl}
In this case $\gg_{\bar{0}} \cong \gs \gl (n) \oplus \gs \gl (n)$.  
By abuse of notation denote by 
$\{ \varepsilon_1, \ldots ,\varepsilon_n, \delta_1, \ldots ,\delta_n\} $
the images of the standard basis elements of $\gh_{\gg \gl (n|n)}^*$
 under the canonical projection $\gh_{\gg \gl (n|n)}^* \mapsto \gh_{\gs \gl (n|n)}^*$.  
Thus,  $\varepsilon_1 + \ldots + \varepsilon_n = 
\delta_1 + \ldots + \delta_n$. Then 
\begin{equation} \label{eqn_psl}
\gh_{\gp \gs \gl (n|n)}^* = \left\{ 
\sum_{i=1}^n a_i \varepsilon_i + \sum_{j=1}^n b_j \delta_j \; | \; \sum_{i=1}^n a_i + \sum_{j=1}^n b_j = 0\right\}. 
\end{equation}
Consider the natural surjective linear map
$\gh_{\gg \gl (n|n)}^* \to \gh_{\gp \gs \gl (n|n)}^*$ defined in terms of 
(\ref{eqn_psl}) by 
\begin{equation}
\label{surject}
\varepsilon_i \mapsto \varepsilon_i - \sigma, 
\; \; 
\delta_j \mapsto \delta_j - \sigma, 
\quad 1 \leq i, j \leq n,
\end{equation}
where $\sigma = \left( \sum_{i=1}^{n} \varepsilon_i +  \sum_{i=1}^{n} \delta_i \right)/(2n)$. 
Then, under the map (\ref{surject}), the set of root of $\gp \gs \gl (n|n)$ is 
identified with the one of $\gs \gl (n|n)$ (and, hence, of $\gg \gl (n|n)$ as well).
Define
\begin{eqnarray*}
L_{\gp \gs \gl (n|n)} (m_0| n_0)& := &  \{ \varepsilon_i - \varepsilon_j \; | \; i, j 
\leq m_0 \mbox{ or } i,j > m_0 \} 
\\
&& \sqcup \{ \delta_k - \delta_l,    
\; | \; k, l \leq n_0 \mbox{ or } k,l > n_0 \} \\
&& \sqcup \{ \pm ( \varepsilon_i - \delta_l)  \; | \; i \leq m_0, l \leq n_0 
\mbox{ or } i >  m_0, l>n_0 \}
\; \; \mbox{and} \\
N^+_{\gp \gs \gl (n|n)} (m_0|n_0) & := & \{ \varepsilon_i - \varepsilon_j \; | 
\; i \leq m_0 <j \} 
\sqcup \{ \delta_k - \delta_l,    \; | \; k \leq n_0<l \}
\\
&& \sqcup \{\varepsilon_i - \delta_l, \delta_k - \varepsilon_j  \; | 
\; i \leq m_0< j; k \leq n_0 < l\}.
\end{eqnarray*}
Set $P_{\gp \gs \gl (n|n)} (m_0, n_0) := L_{\gp \gs \gl (n|n)} (m_0, n_0) 
\sqcup N^+_{\gp \gs \gl (n|n)} (m_0, n_0)$. 
Similarly to the case of $\gg = \gs \gl (m|n)$, 
for $m_0=n_0=0$ and $m_0=n_0=n$, 
$P_{\gp \gs \gl (n|n)} (m_0 |n_0) = \Delta$, and for all other pairs $(m_0, n_0)$ 
$P_{\gp \gs \gl (n|n)} (m_0| n_0)$ is a proper subset of $\Delta$.

Recall that the Lie superalgebras $\gp \gs \gl (n|n)$ are not regular Kac-Moody 
superalgebras. Although there are parabolic subsets of $\gp \gs \gl (n|n)$ 
that are not principal parabolic, every parabolic subset $P$ 
of $\Delta_{\gp \gs \gl (n|n)}$ is the image of 
a parabolic subset $\widehat{P}$ of $\Delta_{\gg \gl (n|n)}$ 
under the map (\ref{surject}),
see \cite[\S 3]{DFG}. Since $\gg \gl (n|n)$ is a quasisimple regular 
Kac-Moody superalgebra, Proposition \ref{kac-moody} applies to $\widehat{P}$. 
Furthermore, using the commutation relations in $\gg \gl (n|n)$, one can verify 
that Proposition \ref{appendix2} is valid for $\gg = \gg \gl (n|n)$, 
$n \geq 2$. Although this proposition fails for $\gg = \gp \gs \gl (3|3)$ 
(take for example: $\alpha = \varepsilon_1 - \delta_1$, 
$\beta =  \varepsilon_2 - \delta_2$, $\alpha + \beta = \delta_3 - \varepsilon_3$), 
we have the following modification.  

\begin{lemma}
\label{psl2-com}
Let $\gg = \gp \gs \gl (3|3)$ and $\alpha, \beta \in \Delta$ be such that 
$\alpha + \beta \neq 0$. Then $[\gg^{\alpha}, \gg^{\beta}] \neq 0$ if and only 
if $\alpha + \beta \in \Delta_{\gg \gl (n|n)}$. Here the root systems of 
$\gh_{\gp \gs \gl (n|n)}^*$ and $\gh_{\gg \gl (n|n)}^*$ are identified via 
the map (\ref{surject}) and the sum $\alpha + \beta$ 
is taken in $\gh_{\gg \gl (n|n)}^*$.
\end{lemma}
\noindent 
{\bf Proof.}  It is sufficient to show that, under the above conditions 
on $\alpha, \beta \in \Delta$, there exist 
$x_{\alpha} \in \gg^{\alpha}$ and $x_{\beta} \in \gg^{\beta}$ such that 
$[x_{\alpha}, x_{\beta}] \neq 0$. We can choose any nonzero elements 
of $\gg \gl (3|3)^{\alpha}$ and $\gg \gl (3|3)^{\beta}$ and consider them 
as elements of $\gg^{\alpha}$ and  $\gg^{\beta}$, respectively. Then using 
the fact that Proposition \ref{appendix2} is valid for $\gg \gl (3|3)$, 
we complete the proof. \hfill $\square$

In view of the above lemma, to obtain a classification of the cominuscule parabolic subsets 
of $\gp \gs \gl (n|n)$, one has to modify the proof of Theorem \ref{prop_slmn} 
for $\gg = \gg \gl (n|n)$ and then transfer the classification 
to $\gg = \gp \gs \gl (n|n)$. The details are left to the reader.

\begin{theorem} \label{prop_psl}
(i) Let $n>2$. Every cominuscule parabolic set of roots of $\gg = \gp \gs \gl (n|n)$ 
is conjugated under the action of the Weyl group 
$W_{\gs \gl (n)} \times W_{\gs \gl (n)}$ of $\gg_{\bar{0}}$ to a unique
subset of the form $P_{\gp \gs \gl (n|n)} (m_0|n_0)$, such that 
$0 \leq m_0 \leq n$, $0 \leq n_0 \leq n$ and $(m_0, n_0) \neq (0,0), (n,n)$. 
Furthermore, for the Levi subalgebras and nilradicals of the corresponding 
parabolic subalgebras we have
$\gl \cong \gs \gl (m_0 | n_0 ) \oplus \gs \gl (n-m_0| n-n_0)$ 
and $\gn^+ \cong V^{m_0|n_0} \otimes \left( V^{n-m_0|n-n_0}\right)^*$
as $\gl$-modules.

(ii) Every cominuscule parabolic set of roots of $\gg = \gp \gs \gl (2|2)$ 
is conjugated under the action of 
the Weyl group $W_{\gs \gl (2)} \times W_{\gs \gl (2)}$ 
of $\gg_{\bar{0}}$ to $P_{\gp \gs \gl (2|2)} (1|1)$ with  
$\gl \cong \gs \gl (1 | 1 ) \oplus \gs \gl (1| 1)$ 
and $\gn^+ \cong V^{1|1} \otimes \left( V^{1|1}\right)^*$
as an $\gl$-module.
\end{theorem}
\subsection{$\gg = \go \gs \gp (2m+1|2n), m \geq 1$.}
\label{osp_odd}
In this case $\gg_{\bar{0}} \cong \gs \go (2m+1) \oplus \gs \gp (2n )$ with
$\Delta_{\bar{0}} = \{ \pm \varepsilon_i \pm \varepsilon_j,  
\pm \varepsilon_i, \pm \delta_k \pm \delta_l, 
\pm  2 \delta_k \; | \; 1 \leq i \neq j \leq m, 1 \leq k \neq l \leq n\}$ 
and $\Delta_{\bar{1}} = \{ \pm \varepsilon_i \pm \delta_k, 
\pm \delta_k \; | \; 1 \leq i \leq m, 1 \leq k \leq n\}$. 
Set
\begin{eqnarray*}
L_{\go \gs \gp (2m+1|2n)} & := &  \{ \pm \varepsilon_i \pm  \varepsilon_j, 
\pm  \varepsilon_i \; | \; j>i>1 \} \sqcup \{ \pm \delta_k \pm \delta_l, 
\pm 2 \delta_k   \; | \; k \neq l\} \\
&& \sqcup
\{ \pm \delta_k, \pm \varepsilon_i  \pm \delta_k,   \; | \; i>1\} \mbox{ and } \\
N^+_{\go \gs \gp (2m+1|2n)} & := & \{ \varepsilon_1 \pm \varepsilon_j, 
\varepsilon_1 \; | \; j>1 \} \sqcup \{ \varepsilon_1 \pm \delta_k  \; | 
\; k  \geq 1\}.
\end{eqnarray*}
and $P_{\go \gs \gp (2m+1|2n)}:= L_{\go \gs \gp (2m+1|2n)} 
\sqcup N^+_{\go \gs \gp (2m+1|2n)}$. 

\begin{theorem} 
\label{prop_osp_odd}
Every cominuscule parabolic set of roots $P$ of $\go \gs \gp (2m+1|2n)$ is conjugated 
under the  action of the Weyl group  $W_{\gs \go (2m+1)} \times W_{\gs \gp (2n)}$ 
of $\gg_{\bar{0}}$ to $ P_{\go \gs \gp (2m+1|2n)}$. For the Levi subalgebra
and nilradical of the corresponding parabolic subalgebra, we have 
$\gl \cong \go \gs \gp (2m -1| 2 n ) \oplus \C$ 
and $\gn^+ \cong V^{2m-1|2n}$ as an $\gl$-module.
\end{theorem}
\noindent 
{\bf Proof.}  
Let $N^+_{\mathfrak{so}}:=N^+ \cap \{ \pm \varepsilon_i, \pm \varepsilon_i 
\pm \varepsilon_j  \; | \; 1 \leq i \neq j \leq m \}$ and 
$N^+_{\mathfrak{sp}}:=N^+ \cap \{ \pm \delta_k \pm \delta_l, 
\pm 2 \delta_k  \; | \; 1 \leq k \neq l \leq n \}$. Obviously, 
$N^+_{\bar{0}} = N^+_{\mathfrak{so}} \sqcup N^+_{\mathfrak{sp}} $.

\medskip
\noindent {\it Case 1: $N^+_{\mathfrak{sp}} \neq \emptyset$}. 
By Proposition \ref{algebras} $P$ is conjugated under the action of the Weyl group 
of $\gs \gp (n)$ to a parabolic subset subset 
with $N^+_{\mathfrak{sp}} = \{ \delta_k + \delta_l, 2\delta_k \; | \; 1 \leq k \neq l \leq n\} $. 
In particular, $\delta_k$ and $2\delta_k$ are in $N^+$ which by Lemma \ref{roots} contradicts 
the assumption that $P$ is cominuscule.

\medskip
\noindent 
{\it Case 2: $N^+_{\mathfrak{s0}}  \neq \emptyset$ 
and $N^+_{\mathfrak{sp}} =  \emptyset$}. Using  Proposition \ref{algebras} again we have that up 
to the action of  the Weyl group of $\gs \go (2m+1)$, 
\begin{eqnarray*}
L_{\bar{0}} & = &  \{\pm \varepsilon_i \pm \varepsilon_j, 
\pm \varepsilon_i\; | \; j> i>1  \sqcup \{ \delta_k - \delta_l \; | \;  
k,l \leq n_0 \mbox { or } k,l>n_0\}  \mbox{ and } \\
N^+_{\bar{0}} & = &  \{ \varepsilon_i - \varepsilon_j, 
\delta_k - \delta_l \; | \; i \leq m_0 < j; \; k \leq n_0 <l\}.
\end{eqnarray*}
Lemma \ref{sums} implies that  $\pm \delta_k$ and $\pm \varepsilon_i  \pm \delta_k$ are in $L$ 
and thus $\varepsilon_1 \pm \delta_k \in N^+$. 
Hence $L = L_{\go \gs \gp (2m+1 | 2n)}$ and  $N^+ = N^+_{\go \gs \gp (2m+1 | 2n)}$.

\medskip
\noindent 
{\it Case 3: $N^+_{\mathfrak{so}} = N^+_{\mathfrak{sp}} =  \emptyset$.} In this case $P$ is not a proper 
subset of $\Delta$.  
\hfill $\square$
\subsection{$\gg = \go \gs \gp (1|2n)$.}
\label{osp1}

In this case $\gg_{\bar{0}} \cong \gs \gp (2n)$ with
$\Delta_{\bar{0}} = \{  \pm \delta_k \pm \delta_l, 
\pm  2 \delta_k \; | \; 1 \leq k \neq l \leq n\}$ 
and $\Delta_{\bar{1}} = \{ \pm \delta_k \; | \; 1 \leq k \leq n\}$.

\begin{theorem} \label{prop_osp_1}
There are no cominuscule parabolic sets of roots of $\go \gs \gp (1|2n)$.
\end{theorem}
\noindent 
{\bf Proof.} If $N^+ \cap \Delta_{\bar{0}} \neq \emptyset$, using the same 
reasoning as in case 1 in the proof of Theorem \ref{prop_osp_odd} 
we reach a contradiction. In the case $N^+ \cap \Delta_{\bar{0}} = \emptyset$ 
one easily proves that $P = L = \Delta$. \hfill $\square$
\subsection{$\gg = \go \gs \gp (2m|2n), m>1$.}
\label{osp_even}
We have $\gg_{\bar{0}} \cong \gs \go (2m) \oplus \gs \gp (2n)$ with
$\Delta_{\bar{0}} = \{ \pm \varepsilon_i  \pm \varepsilon_j, 
\pm \delta_k \pm \delta_l, \pm  2 \delta_k \; | \; 
1 \leq i \neq j \leq m, 1 \leq k \neq l \leq n\}$ 
and $\Delta_{\bar{1}} = \{ \pm \varepsilon_i \pm \delta_k \; | 
\; 1 \leq i \leq m, 1 \leq k \leq n\}$. 
Define the following sets of roots:
\begin{eqnarray*}
L_{\go \gs \gp (2m|2n)}(m)  &:=&   \{  \varepsilon_i -  \varepsilon_j\; | 
\; i \neq j \} \sqcup \{ \delta_k - \delta_l\; | \; k \neq l\} \sqcup 
\{  \pm (\varepsilon_i  - \delta_k)  \; | \; i,k \geq 1\},\\
N^+_{\go \gs \gp (2m|2n)}(m)  &:=&  \{ \varepsilon_i + \varepsilon_j,  \; | 
\; i \neq j \} \sqcup \{\delta_k + \delta_l, 2 \delta_k  \; | \; k \neq l\} 
\sqcup  \{\varepsilon_i +  \delta_k    \; | \; i, k  \geq 1\},\\
L_{\go \gs \gp (2m|2n)}(1)  &:=&   \{  \pm \varepsilon_i \pm  \varepsilon_j\; | 
\; j> i > 1 \} \sqcup \{\pm \delta_k \pm \delta_l, \pm 2 \delta_k\; | \; k \neq l\} 
\\
&& \sqcup \{ \pm \varepsilon_i \pm  \delta_k \; | \; i > 1, k\geq 1\},
\\
N^+_{\go \gs \gp (2m|2n)}(1)  &:=&   \{\varepsilon_1 \pm  
\varepsilon_j    \; | \; j >  1\} \sqcup \{ \varepsilon_1 \pm  \delta_k \; | 
\; k \geq 1\}, \\
\overline{L}_{\go \gs \gp (2m|2n)}(m)  &:=& \bar{\theta} 
L_{\go \gs \gp (2m|2n)}(m), \\
\overline{N}^+_{\go \gs \gp (2m|2n)}(m)  &:=& \bar{\theta} 
N^+_{\go \gs \gp (2m|2n)}(m),  
\end{eqnarray*}
where $\bar{\theta}$ is the involutive automorphism of the Dynkin diagram 
corresponding to the base 
$ \{ \delta_1 - \delta_2,\ldots,\delta_n - \varepsilon_1, 
\varepsilon_1 - \varepsilon_2,\ldots, \varepsilon_{m-1} - \varepsilon_m, 
\varepsilon_{m-1} + \varepsilon_m\}$ that interchanges the last two simple 
roots and preserves all other roots (cf. Proposition \ref{algebras}). 
\begin{theorem} 
\label{prop_osp_even}
There are three orbits of cominuscule parabolic sets of roots of $\go \gs \gp (2m|2n)$ 
under the action of the Weyl group $W_{\gs \go (2m)} \times W_{\gs \gp (2n)}$ of $\gg_{\bar{0}}$. 
The parabolic subsets $P$ with the corresponding Levi superalgebras 
$\gl$ and nilradicals $\gn^+$ (considered as $\gl$-modules) are listed below.

(i) $P_{\go \gs \gp (2m|2n)}(m) := L_{\go \gs \gp (2m|2n)}(m) \sqcup 
N^+_{\go \gs \gp (2m|2n)}(m) $
with $\gl \cong \gg \gl (m| n )$ and $\gn^+ \cong \bigwedge^2 V^{m|n}$.

(ii) $P_{\go \gs \gp (2m|2n)} (1) := L_{\go \gs \gp (2m|2n)}(1)  \sqcup 
N^+_{\go \gs \gp (2m|2n)}(1) $
with $\gl \cong \go \gs  \gp (2m-2| 2n ) \oplus \C$ and $\gn^+ \cong 
V^{2m-2|2n}$.

(iii) $\overline{P}_{\go \gs \gp (2m|2n)}(m)  = 
\overline{L}_{\go \gs \gp (2m|2n)}(m)  \sqcup 
\overline{N}^+_{\go \gs \gp (2m|2n)}(m) $ with  $\gl \cong \gg \gl (m| n )$ 
and $\gn^+ \cong \bigwedge^2 V^{m|n}$. 
\end{theorem}
\noindent 
{\bf Proof.} 
Let $\Delta_{\mathfrak{so}(2m)} =  \{ \pm \varepsilon_i \pm  
\varepsilon_j  \; | \; 1 \leq i \neq j \leq m \}$, $\Delta_{\mathfrak{sp}(2n)} = 
\{\pm  \delta_k \pm \delta_l, 2 \pm \delta_k  \; | \; 1 \leq k \neq l \leq n \}$.  
Define $L_{\mathfrak{so}}:=L \cap \Delta_{\mathfrak{so}(2m)}, N^+_{\mathfrak{so}}:=
N^+ \cap \Delta_{\mathfrak{so}(2m)} $ and $L_{\mathfrak{sp}}:=L \cap \Delta_{\mathfrak{sp}(2n)}, N^+_{\mathfrak{sp}}:=N^+ \cap 
\Delta_{\mathfrak{sp}(2n)}$. Obviously,  $L_{\bar{0}} = L_{\mathfrak{so}} 
\sqcup L_{\mathfrak{sp}} $ and 
$N^+_{\bar{0}} = N^+_{\mathfrak{so}} \sqcup N^+_{\mathfrak{sp}} $.

By Proposition \ref{kac-moody}, we know that $P= P(\Lambda)$ is a principal 
parabolic set of roots determined by some functional $\Lambda$. Let 
$\Lambda (\varepsilon_i ) = x_i,  \Lambda (\delta_k) = y_k$ 
for some $x_i$ and $y_k$.

\medskip
\noindent {\it Case 1: $N^+_{\mathfrak{so}} \neq \emptyset$ and 
$N^+_{\mathfrak{sp}} \neq \emptyset$}.  Proposition \ref{algebras} (iii) 
implies that in the $W_{\gs \gp (2n)}$-orbit of $P$ we have a parabolic subset 
with $L_{\mathfrak{sp}} = \{ \delta_k - \delta_l\; | \; k \neq l\}$ and 
$N^+_{\mathfrak{sp}} = \{ \delta_k + \delta_l, 2 \delta_k \; | \; k \neq l\}$.
Next,  following Proposition \ref{algebras} (iv), we consider three sub-cases 
for the $W_{\gs \go (2m)}$-orbit of $P$.

\medskip
\noindent {\it Case 1.1: $L_{\mathfrak{so}} = L_{\gs \go (2m)}(1)$ and 
$N^+_{\mathfrak{so}} = N_{\gs \go (2m)}^+(1)$.}  In this case $x_2 =  \ldots = x_m = 0$, 
$x_1 >0$, and $y_1 = \ldots =y_n >0$. But then $\varepsilon_2 + \delta_1$ and 
$-\varepsilon_2 + \delta_1$ are in $N^+$ and 
$(\varepsilon_2 + \delta_1) + (-\varepsilon_2 + \delta_1) \in \Delta$, 
which contradicts the assumption that $P$ is cominuscule. 

\medskip
\noindent {\it Case 1.2: $L_{\mathfrak{so}} = L_{\gs \go (2m)}(m)$ and 
$N^+_{\mathfrak{so}} = N^+_{\gs \go (2m)}(m)$.}  Now  
$x_1 = x_2 = \ldots = x_m = x > 0$  and $y_1 = \ldots =y_n = y >0$. If $x>y$ then 
$\varepsilon_1 - \delta_1, 2 \delta_1 \in N^+$ with  
$(\varepsilon_1 - \delta_1) + 2 \delta_1 \in \Delta$, which contradict 
the assumption that $P$ is cominuscule. Similarly if $y>x$ we use that 
$\delta_1 \pm  \varepsilon_1 \in N^+$ and reach a contradiction. It remains 
to consider $x = y$. In this case one verifies that 
$L = L_{\go \gs \gp (2m|2n)}(m) $ and $N^+ = N^+_{\go \gs \gp (2m|2n)}(m)$.

\medskip
\noindent {\it Case 1.3: $L_{\mathfrak{so}}  = \bar{\theta} L_{\gs \go (2m)}(m)$ 
and $N^+_{\mathfrak{so}}  = \bar{\theta} N^+_{\gs \go (2m)}(m)$. } This follows 
from Proposition \ref{algebras} and Case 1.2.  We obtain 
$L = \bar{\theta} L_{\go \gs \gp (2m|2n)}(m)$ 
and $N^+ = \bar{\theta} N^+_{\go \gs \gp (2m|2n)}(m)$.

\medskip
\noindent {\it Case 2: $N^+_{\mathfrak{so}}  \neq \emptyset$ and 
$N^+_{\mathfrak{sp}} = \emptyset$}.  In this case  $y_1 = \ldots =y_n =0$. 
We consider again three subcases for $N^+_{\mathfrak{so}} $.

\medskip
\noindent {\it Case 2.1: $L_{\mathfrak{so}}  = L_{\gs \go (2m)}(1)$ and 
$N^+_{\mathfrak{so}} = N^+_{\gs \go (2m)}(1)$}. We have $x_2 = \ldots = x_m =0$ 
and $x_1>0$. This leads to $L = L_{\go \gs \gp (2m|2n)}(1) $ and 
$N^+ = N^+_{\go \gs \gp (2m|2n)}(1)$.

\medskip
\noindent {\it Case 2.2: $L_{\mathfrak{so}}  = L_{\gs \go (2m)}(m)$ and 
$N^+_{\mathfrak{so}} = N_{\gs \go (2m)}^+(m)$}. In this case 
$x_1 = \ldots = x_m > 0$. In this case 
$\varepsilon_1 - \delta_1, \varepsilon_2 + \delta_1 \in N^+$ with 
$(\varepsilon_1 - \delta_1)+(\varepsilon_2 + \delta_1) \in \Delta$, 
which leads to a contradiction.

\medskip
\noindent {\it Case 2.3: $L_{\mathfrak{so}}  = \bar{\theta} L_{\gs \go (2m)}(m)$ 
and $N^+_{\mathfrak{so}} = \bar{\theta} N^+_{\gs \go (2m)}(m)$}. Using Case 2.2 
and Proposition \ref{algebras}(iv) we reach again a contradiction.

\medskip
\noindent {\it Case 3: $N^+_{\mathfrak{so}}  = \emptyset$ and 
$N^+_{\mathfrak{sp}}  \neq \emptyset$}.  In this case $x_1= \ldots =x_m = 0$ and 
$y_1= \ldots =y_n>0$. We now have that 
$\varepsilon_1 + \delta_1, - \varepsilon_1 + \delta_1 \in N^+$ 
with $(\varepsilon_1 + \delta_1)+( - \varepsilon_1 + \delta_1) \in \Delta$, 
again  a contradiction.

\medskip
\noindent {\it Case 4: $N^+_{\mathfrak{so}} =  N^+_{\mathfrak{sp}} = \emptyset$}. 
In this case $P = L = \Delta$.  \hfill $\square$
\subsection{$\gg = \go \gs \gp (2|2n)$.}
\label{osp2}
We have $\gg_{\bar{0}} \cong \C \oplus \gs \gp (2n)$ with
$\Delta_{\bar{0}} = \{ \pm \delta_k \pm \delta_l, 
\pm  2 \delta_k \; | \;1 \leq k \neq l \leq n\}$ and $\Delta_{\bar{1}} = 
\{ \pm \varepsilon_1 \pm \delta_k \; | \; 1 \leq k \leq n\}$. 
Let us define the following sets of roots:
\begin{eqnarray*}
L_{\go \gs \gp (2|2n)} (0) & := &  \Delta_{\bar{0}}, \\
N^+_{\go \gs \gp (2|2n)}(0) & := & \Delta_{\bar{1}}^{+} = 
\{ \varepsilon_1 \pm \delta_k    \; | \; k  \geq 1\}, \\
L_{\go \gs \gp (2|2n)} (n) & := & \{\delta_k - \delta_l \; | 
\; k \neq l\} \sqcup \{\pm (\varepsilon_1 - \delta_k),   \; | 
\; k\geq 1\},\\
N^+_{\go \gs \gp (2|2n)}(n) & := & \{ \delta_k + \delta_l, 2 \delta_k \; | 
\; k \neq l \} \sqcup \{ \varepsilon_1 + \delta_k    \; | \; k  \geq 1\},\\
\overline{L}_{\go \gs \gp (2|2n)} (n) & := & \{\delta_k - \delta_l \; | 
\; k \neq l\} \sqcup \{\pm (\varepsilon_1 + \delta_k),   \; | 
\; k\geq 1\},\\
\overline{N}^+_{\go \gs \gp (2|2n)}(n) & := & \{ \delta_k + \delta_l, 
2 \delta_k \; | \; k \neq l \} \sqcup \{ - \varepsilon_1 + \delta_k    \; | 
\; k  \geq 1\}.
\end{eqnarray*}
Note that $\overline{L}_{\go \gs \gp (2|2n)} (n)  = 
\tau L_{\go \gs \gp (2|2n)} (n)$ and 
$\overline{N}^+_{\go \gs \gp (2|2n)} (n)  = \tau N^+_{\go \gs \gp (2|2n)} (n)$, 
where $\tau$ is the automorphism of the Dynkin diagram of the base 
$\{ \delta_1 - \delta_2, \ldots ,\delta_{n-1} - \delta_{n}, \delta_{n} - \varepsilon_1, 
\delta_{n} + \varepsilon_1\}$ of $\Delta$ that interchanges the last two roots 
and preserves all other roots (one can think of $\tau$ as the diagram automorphism 
that interchanges $\varepsilon_1$ with $- \varepsilon_1$). Note also that $\tau \Delta_{\bar{1}}^{+} = - \Delta_{\bar{1}}^{+}$. 

Set $P_{\go \gs \gp (2|2n)}(0) := L_{\go \gs \gp (2|2n)} (0) \sqcup 
N^+_{\go \gs \gp (2|2n)}(0)$,  $P_{\go \gs \gp (2|2n)}(n) := 
L_{\go \gs \gp (2|2n)} (0) \sqcup N^+_{\go \gs \gp (2|2n)}(n)$, 
and $\overline{P}_{\go \gs \gp (2|2n)}(n) := 
\overline{L}_{\go \gs \gp (2|2n)} (n) \sqcup 
\overline{N}^+_{\go \gs \gp (2|2n)}(n)$.

\begin{theorem} 
\label{prop_osp2}
There are  four orbits of cominuscule parabolic sets of roots 
of $\gg = \go \gs \gp (2|2n)$ under the action of the  Weyl group 
$W_{\gs \gp (2n)}$ of $\gg_{\bar{0}}$: $P_{\go \gs \gp (2|2n)}(0)$, 
$-P_{\go \gs \gp (2|2n)}(0)$,\footnote{one can certainly consider 
$- P_{\go \gs \gp (2|2n)}(0)$ as $\tau P _{\go \gs \gp (2|2n)}(0)$.}  
$P_{\go \gs \gp (2|2n)}(n)$,  and $\overline{P}_{\go \gs \gp (2|2n)}(n)$. 
The corresponding Levi subalgebras $\gl$ and nilradicals $\gn^+$ 
(considered as $\gl$-modules) are given by  
$\gl \cong \gs \gp (2 n ) \oplus \C, \gn^+ \cong V^{2n}$ 
in the first two cases, and 
$\gl \cong \gs \gl (1| n ), \gn^+ \cong S^2 V^{1|n}$ 
in the last two. 
\end{theorem}
\noindent 
{\bf Proof.} 
By Proposition \ref{kac-moody} we have $P= P(\Lambda)$ 
for some functional $\Lambda$. Let $\Lambda (\varepsilon_1) = x,  
\Lambda (\delta_k) = y_k$, $k=1, \ldots ,n$, for some $x$ and $y_k$ in $\C$.

\medskip
\noindent {\it Case 1: $N^+_{\bar{0}} \neq \emptyset$}. Due to 
Proposition \ref{algebras}  (iii) we may assume that 
$L_{\bar{0}} = L_{\gs \gp (2n)}$ and $N^+_{\bar{0}} = N^+_{\gs \gp (2n)}$.  
Thus $y_1 =  \ldots = y_n = y>0$. We proceed with three sub-cases.

\medskip
\noindent {\it Case 1.1: $y>|x|$}. In this case $\pm \varepsilon + \delta_1$ 
are in $N^+$ and their sum is a root, which is a contradiction. 

\medskip
\noindent {\it Case 1.2: $y<|x|$}. Similarly to the previous case we 
reach a contradiction with the assumption that $P$ is cominuscule. 
For example, if  $-x>y>x>0$ we have that $-\delta_1 - \varepsilon_1$ 
and $2 \delta_1$ are in $N^+$. The other cases are analogous. 

\medskip
\noindent {\it Case 1.3: $y=|x|$}. It easily follows that in the case $y=x$ 
we obtain $P = P_{\go \gs \gp (2|2n)}(n)$, while in the case $y=-x$ we find 
$P = \overline{P}_{\go \gs \gp (2|2n)}(n)$. 

\medskip
\noindent {\it Case 2: $N^+_{\bar{0}} = \emptyset$}. In this case 
$y_1= \ldots =y_n=0$, 
and in particular $L = \Delta_{\bar{0}}$. Depending on the sign of $x$ 
we have $P = P_{\go \gs \gp (2|2n)}(0) $ 
or $P = -P_{\go \gs \gp (2|2n)}(0)$. \hfill $\square$

\subsection{$\gg = D(2,1; \alpha)$.}
\label{D_alpha}
We have $\gg_{\bar{0}} \cong \gs \gl (2) \oplus \gs \gl (2) \oplus \gs \gl (2)$ 
with $\Delta_{\bar{0}} = \{   \pm \gamma_i \; | \;i = 1,2,3\}$ and 
$\Delta_{\bar{1}} = \{  \frac{1}{2}( \pm \gamma_1 \pm \gamma_2 \pm \gamma_3)\}$.  
Here $\pm \gamma_i$ denote the roots of the $i$-th copy of $\gs \gl (2)$ in 
$\gg_{\bar{0}} $.  In this case we will classify the cominuscule parabolics of $\gg$ 
up to the action of the group $S(\Delta) = W_{\gg_{\bar{0}}} \ltimes S_3$, 
where $S_3$ acts by permutations on $\{ \gamma_1, \gamma_2, \gamma_3\}$. 
(In fact $S(\Delta)$ is the group of automorphisms of $\Delta$.) 
By considering the action of a larger group than $W_{\gg_{\bar{0}}}$ we 
avoid a longer list of similarly behaved cominuscule parabolics. We leave to the 
reader to reconstruct from this the $W_{\gg_{\bar{0}}}$-orbits of cominuscule 
parabolic subsets of roots of $D(2,1; \alpha)$.

The root system of $\gg$ coincides with the root system of $\go \gs \gp (4|2)$ 
and an explicit isomorphism $\Delta_{\gg} \to \Delta_{\go \gs \gp (4|2)}$ 
is provided by: $\gamma_1 \mapsto \varepsilon_1 +  \varepsilon_2$, 
$\gamma_2 \mapsto \varepsilon_1 -  \varepsilon_2$, $\gamma_3 \mapsto 2 \delta_1$. 
Using this equivalence and Theorem \ref{prop_osp_even} 
one easily verifies the following.

\begin{theorem} 
\label{prop_d_alpha}
There is only one orbit of cominuscule parabolic sets of roots 
of $D(1,2; \alpha)$ 
under the action of  $S(\Delta)$: 
$P_{D(1,2; \alpha)} =L_{D(1,2; \alpha)} \sqcup N_{D(1,2; \alpha)} $, 
where $L_{D(1,2; \alpha)} = \{ \pm \gamma_3, \frac{1}{2}(\pm(\gamma_1 - \gamma_2) 
\pm \gamma_3)\}$ and $N^+_{D(1,2; \alpha)}  = \{ \gamma_1, \gamma_2, 
\frac{1}{2}(\gamma_1 + \gamma_2 \pm \gamma_3)\}$. For the corresponding Levi subalgebra 
and nilradical we have $\gl \cong \gg \gl (2|1 )$ and 
$\gn^+ \cong \bigwedge^2 V^{2|1}$ as an $\gl$-module.
\end{theorem}

\subsection{$\gg = F(4)$.}
\label{F4}
We have $\gg_{\bar{0}} \cong \gs \go (7) \oplus \gs \gl (2)$ 
with $\Delta_{\bar{0}} = \{   \pm \varepsilon_i \pm \varepsilon_j,  \pm \varepsilon_i, \pm \gamma \; | \; 1\leq i \neq j \leq 3 \}$ and 
$\Delta_{\bar{1}} = \left\{  \frac{1}{2}\left( \pm \varepsilon_1 \pm \varepsilon_2 \pm \varepsilon_3   \pm \gamma \right)\right\}$.  
Here $\{  \pm \varepsilon_i \pm \varepsilon_j,  \pm \varepsilon_i \; | \; 1\leq i \neq j \leq 3 \}$ denote the roots of $\gs \go (7)$, while $\pm \gamma$ denote the roots of $\gs \gl (2)$ in 
$\gg_{\bar{0}} $.  

\begin{theorem} 
\label{prop_f4}
There are no  cominuscule parabolic sets of roots of $F(4)$.
\end{theorem}
\noindent 
{\bf Proof.} 
By Proposition \ref{kac-moody} we have $P= P(\Lambda)$ 
for some functional $\Lambda$. Let $\Lambda (\varepsilon_i) = x_i,  
\Lambda (\gamma) = y$, $i=1, \ldots ,3$, for some $x_i$ and $y$ in $\C$.

\medskip
\noindent {\it Case 1: }$N^+ \cap \{   \pm \varepsilon_i \pm \varepsilon_j,  \pm \varepsilon_i, \pm \gamma \; | \; 1\leq i \neq j \leq 3 \} \neq \emptyset$.  Due to 
Proposition \ref{algebras}  (ii) we may assume that $N^+ \cap \{   \pm \varepsilon_i \pm \varepsilon_j,  \pm \varepsilon_i, \pm \gamma \; | \; 1\leq i \neq j \leq 3 \}  = N^+_{\gs \go (5)}$. Thus $x_1>0, x_2=x_3=0$. But then $\frac{1}{2} (\varepsilon_1\pm\varepsilon_2 )\in N^+$ and their sum is a root, which implies that $P$ is not cominuscule.

\medskip
\noindent {\it Case 2: }$N^+ \cap \{   \pm \varepsilon_i \pm \varepsilon_j,  \pm \varepsilon_i, \pm \gamma \; | \; 1\leq i \neq j \leq 3 \} = \emptyset$. In this case $x_1=x_2=x_3 = 0$. If $y>0$ then $\frac{1}{2} (\varepsilon_1\pm \gamma ) \in N^+$ which contradicts again to the fact that $P$ is cominuscule. The case $y<0$ is similar to the case $y>0$, while for $y=0$ we obtain $P=\Delta$.
 \hfill $\square$

\subsection{$\gg = G(3)$.}
\label{G3}
We have $\gg_{\bar{0}} \cong G_2 \oplus \gs \gl (2)$ 
with $\Delta_{\bar{0}} = \{   \varepsilon_i - \varepsilon_j,  \pm \varepsilon_i, \pm \gamma \; | \; 1\leq i \neq j \leq 3 \}$ and 
$\Delta_{\bar{1}} = \left\{ \pm  \frac{\gamma}{2}, \pm \varepsilon_i  \pm \frac{\gamma}{2} \; | \; 1\leq i\leq 3  \right\}$.  
Here $\varepsilon_1 + \varepsilon_2 + \varepsilon_3 = 0$ and $\{   \varepsilon_i - \varepsilon_j,  \pm \varepsilon_i \; | \; 1\leq i \neq j \leq 3 \}$ denote the roots of $G_2$, while $\pm \gamma$ denote the roots of $\gs \gl (2)$ in 
$\gg_{\bar{0}} $.  

\begin{theorem} 
\label{prop_g3}
There are no  cominuscule parabolic sets of roots of $G(3)$.
\end{theorem}
\noindent 
{\bf Proof.} 
By Proposition \ref{kac-moody} we have $P= P(\Lambda)$ 
for some functional $\Lambda$. Let $\Lambda (\varepsilon_i) = x_i,  
\Lambda (\gamma) = y$, $i=1, 2$, for some $x_i$ and $y$ in $\C$. Since $G_2$ has no cominuscule parabolic sets of roots we have that $N^+ \cap \{   \varepsilon_i - \varepsilon_j,  \pm \varepsilon_i \; | \; 1\leq i \neq j \leq 3 \}  = \emptyset$. This implies $x_1 = x_2 = 0$. In the case $y >0$ we have that $\pm \varepsilon_1 + \frac{\gamma}{2} \in N^+$ and their sum is a root, which contradicts with the fact that $P$ is cominuscule. The case $y<0$ is similar, while for $y=0$ we obtain $P=\Delta$.
 \hfill $\square$ 
 

\section{Classification of 
the cominuscule parabolics of strange classical Lie superalgebras}
\label{strange}
\setcounter{equation}{0}

In this section we classify the cominuscule parabolics of the two 
strange classical Lie superalgebras $\gg = \gp \gs \gq (n)$ 
and $\gg = {\bf s}\gp (n)$. As in the case of the basic classical
Lie superalgebras, for these superalgebras the even part 
$\gg_{\bar{0}}$ of $\gg$ is a reductive Lie algebra. Analogously to 
the previous section, the Weyl group $W_{\gg_{\bar{0}}}$ of 
$\gg_{\bar{0}}$ acts on the root system $\Delta$ of $\gg$ 
(and thus on $\Delta \cup (- \Delta)$). This action
induces an action of $W_{\gg_{\bar{0}}}$ on the set of parabolic subsets of $\Delta$,
which preserves the class of cominuscule parabolic subsets of $\Delta$. 
We will classify the cominuscule parabolic subsets of $\Delta$ up to this 
action of the Weyl group $W_{\gg_{\bar{0}}}$.

In all proofs we will assume 
that $P = \widetilde{P} \cap \Delta$ is a cominuscule parabolic set of roots 
for some parabolic subset $\widetilde{P}$ of $\Delta \cup (- \Delta)$ 
for which the corresponding nilradical $\gn^+$ of $\gp_P$ is abelian,
recall Definition \ref{dLN}.
The root system of $\gp \gs \gq (n)$ is symmetric 
(and $P = \widetilde{P}$ for $\gg = \gp \gs \gq (n)$), 
while the one of ${\bf s}\gp (n)$ is not. 
We will use the notation  $L$, $N^+$, 
$\gl$, and $\gn^+$ from Definitions \ref{dLN} and \ref{LN-alg}. 
Set $P_{\bar{0}}:=P \cap \Delta_{\bar{0}}$. By Lemma \ref{restrict}, 
either $P_{\bar{0}} = \Delta_{\bar{0}}$ or $P_{\bar{0}}$ 
is a parabolic subset of $\gg_{\bar{0}}$ with Levi component 
$L_{\bar{0}}:=L \cap \Delta_{\bar{0}}$ 
and nilradical $N_{\bar{0}}^+:=N^+ \cap \Delta_{\bar{0}}$. 
We also set $P^- = L \sqcup N^-= (- \widetilde{P})\cap \Delta $
and $N^- = \widetilde{N}^{-} \cap \Delta$
where $\widetilde{N}^{-} = (- \widetilde{P}) \backslash
\widetilde{P}$, cf. Lemma \ref{sums}. Then 
$P^-$ is a parabolic subset of $\Delta$ and 
$P^- = L \sqcup N^-$ is a Levi decomposition of $P^-$.
The corresponding Levi decomposition of the parabolic
subalgebra $\gp_{P^-}$ of $\gg$ is 
$\gp_{P^-} = \gl \ltimes \gn^-$, where
$\gn^-:= \bigoplus_{\alpha \in N^-} \gg^{\alpha}$.
\subsection{$\gg = \gp \gs \gq (n)$.}
\label{q_n}
In this case $\gg_{\bar{0}} \cong \gs \gl (n)$ with
$\Delta_{\bar{0}} = \Delta_{\bar{1}} = \{   \varepsilon_i -  \varepsilon_j \; | 
\;1 \leq i \neq j \leq n\}$. Because $\Delta$ coincides with the root system of 
$\gs \gl_n$ one can easily modify Proposition \ref{algebras} (i) and obtain 
the classification of the cominuscule parabolic subalgebras of $\gp \gs \gq (n)$.  
For $1 \leq n_0 \leq n-1$, set $L_{\gp \gs \gq (n)} (n_0)= L_{\gs \gl (n)} (n_0)$,  
$N^+_{\gp \gs \gq (n)} (n_0)= N_{\gs \gl (n)} (n_0)$, and 
$P_{\gp \gs \gq (n)} (n_0)= P_{\gs \gl (n)} (n_0)$. The details are left to 
the reader. For the next theorem we introduce some notation.  Recall that $\gq(m)$ is the Lie superalgebra of all matrices $X =  \left(  \begin{matrix}A & B \\ B & A  \end{matrix} \right)$ where $A$ and $B$ are $m \times m$ matrices. We set $\mbox{otr} (X):= \tr B$. For $n_0$, $1\leq n_0 \leq n-1$, denote 
$$
\gs \gq (n_0, n-n_0) := \left\{ (X,Y) \in \gq (n_0) \oplus \gq (n-n_0) \; | \; \mbox{otr} (X) + \mbox{otr}(Y) = 0\right\},
$$
and set $\gp \gs \gq (n_0, n-n_0):= \gs \gq (n_0, n-n_0)/ (\C \Id)$. In particular, $\gp \gs \gq (1, n-1) \simeq \gp \gs \gq (n-1, 1) \simeq \gp \gs \gq (n-1)$.
\begin{theorem}
\label{prop_qn} 
There are $n-1$ orbits of cominuscule parabolic 
sets of roots of $\gg =  \gp \gs \gq (n)$ under the action of 
the Weyl group $W_{\gs \gl (n)}$ of $\gg_{\bar{0}}$ with representatives 
$P_{\gp \gs \gq (n)} (n_0)$, for $1 \leq n_0 \leq n-1$. The corresponding
Levi subalgebras $\gl$ and nilradicals (considered as $\gl$-modules) are 
given by $\gl \cong \gp \gs \gq (n_0, n-n_0)$
and $\gn^+ \cong V^{n_0|n_0} \otimes (V^{n-n_0 | n-n_0})^*$. 
\end{theorem}
\subsection{$\gg = {\bf s} \gp (n)$.}
\label{p_n}
We have $\gg_{\bar{0}} \cong \gs \gl (n)$ with
$\Delta_{\bar{0}} = \{   \varepsilon_i -  \varepsilon_j \; | 
\;1 \leq i \neq j \leq n\}$ and $\Delta_{\bar{1}} = 
\{ \pm (\varepsilon_i + \varepsilon_j), 2 \varepsilon_i \; | 
\; 1 \leq i < j  \leq n\}$. Since $\Delta \neq - \Delta$, we need to
pass to the parabolic subsets of $\Delta \cup -\Delta$. On the other 
hand $\Delta \cup (- \Delta)$ coincides with the root system of $\gs \gp (2n)$, 
hence all parabolic subsets of $\Delta \cup (- \Delta)$ are principal,
cf. \S \ref{principal}.
Let us define the following sets of roots:
\begin{eqnarray*}
L_{{\bf s}\gp (n)} (0) & := &  \Delta_{\bar{0}}, \\
N_{{\bf s}\gp (n)}^{+}  (0)& := &  \Delta_{\bar{1}}^{+} =  
\{ \varepsilon_i + \varepsilon_j, 2 \varepsilon_i \; | \; i < j\}, \\
N_{{\bf s}\gp (n)}^{-} (0) & := &  \Delta_{\bar{1}}^{-} =  
\{ - \varepsilon_i - \varepsilon_j \; | \; i < j\}, \\
L_{{\bf s}\gp (n)} (n_0) & := &   \{   \varepsilon_i -  \varepsilon_j \; | 
\;i, j \leq n_0\mbox{ or } i,j>n_0\} \sqcup  
\{ \pm (\varepsilon_i +  \varepsilon_j) \; | \;i \leq n_0 < j\}, \\
N_{{\bf s}\gp (n)}^{+}  (n_0) & := &   \{   \varepsilon_i  -  \varepsilon_j \; | 
\; i \leq n_0 < j \} \sqcup \{   \varepsilon_i  +  \varepsilon_j \; | 
\; i \leq j \leq n_0 \} 
\\
&& \sqcup  \{  -  \varepsilon_i  -  \varepsilon_j \; | 
\;  i >  j > n_0 \}, \\
L_{{\bf s}\gp (n)} (n) & := &   \{   \varepsilon_i -  \varepsilon_j \; | 
\;i, j < n\} \sqcup  
\{ \pm (\varepsilon_i +  \varepsilon_j), 2 \varepsilon_i \; | \;i <j<n\}, \\
N_{{\bf s}\gp (n)}^{+}  (n)  & := &   \{   - \varepsilon_n 
\pm  \varepsilon_j \; | \; j <n\},  
\end{eqnarray*}
for $1 \leq n_0 \leq n-1$. Set 
$P_{{\bf s}\gp (n)} (n_0)  := L_{{\bf s}\gp (n)}(n_0) 
\sqcup N_{{\bf s}\gp (n)} ^{+} (n_0)$ for $0 \leq n_0 \leq n$.

\begin{theorem} 
\label{prop_pn}
There are $n+2$ orbits of cominuscule parabolic sets of roots of 
$\gg = {\bf s} \gp (n)$ under the action of 
the Weyl group $W_{\gs \gl (n)}$ of $\gg_{\bar{0}}$
with representatives 
$P_{{\bf s}\gp (n)} (n_0)$ for $0 \leq n_0 \leq n$ and 
$P_{{\bf s}\gp (n)}^{-} (0) = -P_{{\bf s}\gp (n)}(0)$. 
All cominuscule parabolic sets of roots have unique 
Levi decompositions and the Levi subalgebras $\gl$ and nilradicals
(considered as $\gl$-modules) of the above parabolics are given by:
\begin{itemize}
\item  $\gl \cong \gs \gl (n) , \gn^+ \cong S^2 V^{n}, \gn^{-} 
\cong \bigwedge^2 (V^n)^*$ for $n_0=0$;
\item $\gl \cong \gs \gl (n_0 | n-n_0) , \gn^+ 
\cong S^2 V^{n_0 | n-n_0}$,
for $1 \leq n_0 \leq n-1$;
\item $\gl \cong {\bf s} \gp (n-1 ) \oplus \C , \gn^+ \cong V^{n-1|n-1}$, 
for $n_0=n$.
\end{itemize}
\end{theorem}
\noindent 
{\bf Proof.} 
Since $\Delta \cup (-\Delta)$ coincides with the root system
of $\gs \gp ( 2n)$, all parabolic sets of roots of $\Delta \cup (-\Delta)$
are principal. We consider $\Delta \cup (- \Delta)$ as 
a subset of a real vector space $V$ 
with basis $\{ \varepsilon_i \; | \; 1 \leq i \leq n \}$. We have 
$\widetilde{P}= \widetilde{P}(\Lambda)$ for some functional $\Lambda$ on $V$. 
Let $\Lambda (\varepsilon_i) = x_i$, $1 \leq i \leq n$, for some $x_i \in \R$.  

\medskip
\noindent {\it Case 1: $N_{\bar{0}}^+ \neq \emptyset$}. By 
Proposition \ref{algebras} we have that, up to the action of 
the Weyl group $W_{\gs \gl (n)}$ of $\gg_{\bar{0}}$, 
$L_{\bar{0}} = L_{\gs \gl (n)}(n_0)$ and $N_{\bar{0}}^+ = N_{\gs \gl (n)}(n_0)$ 
for some $n_0$, $1 \leq n_0 \leq n-1$. Thus $x_1= \ldots =x_{n_0} = x$, 
$x_{n_0+1} = \ldots =x_n = y$, for some $x$ and $y$ such that $x>y$. 
We proceed with three subcases.

\medskip
\noindent {\it Case 1.1: $x+y>0$}. We  have that $\varepsilon_1 - \varepsilon_n$ 
and $\varepsilon_1 + \varepsilon_n$ are in $N^+$ and their sum is a root, 
which is a contradiction. 

\medskip
\noindent {\it Case 1.2: $x+y=0$}. 
We have that $x = -y >0$. In this case $L =  L_{{\bf s}\gp (n)}(n_0)$ and 
$N^+ =  N_{{\bf s}\gp (n)}(n_0)$. 

\medskip
\noindent {\it Case 1.3: $x+y<0$}. If $n_0<n-1$, then 
$\varepsilon_1 - \varepsilon_{n-1}$ and $-\varepsilon_{1} - \varepsilon_{n}$ 
are in $N^+$ and their sum is a root. It remains to consider 
the case $n_0 = n-1$. We easily 
see that $\pm (\varepsilon_i + \varepsilon_j)$ are not in $N^+$ for every $i,j>1$. 
Thus $x = 0$ which leads to $L =  L_{{\bf s}\gp (n)}(n)$ and 
$N^+ =  N^+_{{\bf s}\gp (n)}(n)$.

\medskip
\noindent {\it Case 2: $N_{\bar{0}}^+ = \emptyset$}. In this case 
$\Delta_{\bar{0}} \subset L_{\bar{0}}$ and in particular $x_1 = \ldots =x_n= x$. 
If $x>0$ we obtain $P = P_{{\bf s}\gp (n)} (0)$, while for $x<0$ we have 
$P = P^{-}_{{\bf s}\gp (n)} (0)$. 

The isomorphisms for the Levi components $\gl$ and the structure of 
the nilradicals as $\gl$-modules are straightforward and are left to 
the reader.
\hfill $\square$
\\
\noindent
We extend Remark \ref{nonu1} with an illustration of the nonuniquenes 
of Levi decompositions of parabolic sets of roots for the root system 
of ${\bf s} \gp (n)$.
\begin{remark}
\label{nonu2} Denote the subset of roots 
\[
P:= \{ \varepsilon_i \pm \varepsilon_j \; | \; 1 \leq i < j \leq n\} 
\sqcup \{ 2 \varepsilon_i \; | \; 1 \leq i \leq n\}
\subset \Delta. 
\]
Let us identify $V^*$ with $\R^n$, where 
$\Lambda \in V^* \mapsto ( \Lambda (\varepsilon_1), \ldots, \Lambda (\varepsilon_n))$.
The polyhedron ${\mathcal{F}}(P)$ defined in Remark \ref{nonu1} 
is given by 
\[
{\mathcal{F}}( P ) = \{ (x_1, \ldots, x_n) \in \R^n \; | 
\; x_1 > \ldots > x_n \geq 0 \}.
\] 
In particular, ${\mathcal{F}}(P)$ is nonempty and $P$ is 
a principal parabolic subset of $\Delta$.
Furthermore, ${\mathcal{F}}(P)$ has two faces with interiors
\[
\{ (x_1, \ldots, x_n) \in \R^n \; | \; x_1 > \ldots > x_n > 0 \} \; \; 
\mbox{and} \; \; 
\{ (x_1, \ldots, x_n) \in \R^n \; | \; x_1 > \ldots > x_n = 0 \}.
\]
The corresponding Levi components are given by 
\[
L= \emptyset \; \; \mbox{and} \; \; 
L= \{ 2 \varepsilon_n \},
\]
respectively. One easily generalizes this example to show that 
for every parabolic subset $P$ of the root system of 
${\bf s} \gp (n)$, the polyhedron 
${\mathcal{F}}(P)$ has at most two faces. (The small number of faces 
is due to the fact that $\Delta \backslash (- \Delta)$ has only 
$n$ roots.) The corresponding Levi 
decompositions of $P$ are easily described in a similar fashion.
\end{remark} 
\section{Classification of  cominuscule parabolics of Cartan type Lie superalgebras}
\label{Cartan}
\setcounter{equation}{0}

In this section we classify the cominuscule parabolic sets of roots 
of all Cartan type Lie superalgebras $\gg$. Each such Lie superalgebra 
has a natural Cartan subalgebra $\gh$. The corresponding root system 
will be denoted by $\Delta$. In each case the even part $\gg_{\bar{0}}$ 
has a natural Levi subalgebra $\gl_{\bar{0}}$ such that 
$\gl_{\bar{0}} \cap \gh = \gh_{\bar{0}}$. The action of the Weyl group 
$W_{\gl_{\bar{0}}}$ of the Levi subalgebra 
$\gl_{\bar{0}}$ on the root system of $\gl_{\bar{0}}$
extends  to an action of $W_{\gl_{\bar{0}}}$ on $\Delta$. 
For each parabolic subset $P$ of $\Delta$ with 
a Levi decomposition $P = L \sqcup N^+$ and $w \in W_{\gl_{\bar{0}}}$, 
$w(P)$ is a parabolic subset of $\Delta$ and 
$w(P) = w(L) \sqcup w(N^+)$ is a Levi decomposition 
of $w(P)$. Furthermore, $P$ is cominuscule if and only 
if $w(P)$ is cominuscule. Our classification amounts to classifying
the orbits of cominuscule parabolic subsets of $\Delta$ under the
action of the Weyl group $W_{\gl_{\bar{0}}}$ of the Levi subalgebra
$\gl_{\bar{0}}$ of $\gg_{\bar{0}}$.

Given a pair of integers $j \leq k$, set $[j,k]:= \{j,\ldots ,k \}$. 
\subsection{$\gg = W(n)$.}
\label{W_n} 
Let $\gg:=W(n)= W(\xi_1, \ldots,\xi_n)$ be the Lie superalgebra consisting of the superderivations of the Grassmann algebra 
$\bigwedge(n) =\bigwedge(\xi_1, \ldots ,\xi_n)$. 
The elements of $W(n)$ have the 
form $\sum_{i=1}^n p_i(\xi_1, \ldots,\xi_n) \frac{\partial}{\partial \xi_i}$, 
where $p_i \in \bigwedge(n)$ and $\frac{\partial}{\partial \xi_i}$
are the derivations of $\bigwedge(n)$ such that
$\frac{\partial}{\partial \xi_l} (\xi_l) = \delta_{il}$,
for all $1 \leq l \leq n$. 
The standard Cartan subalgebra of $W(n)$ is
\[
\gh = \Span \left\{ \xi_i \frac{\partial}{\partial \xi_i} \; \Big{|} \; 
1 \leq i \leq n \right\}.
\]
Both $\bigwedge(n)$ 
and $W(n)$ have natural gradings: $\bigwedge(n) = \bigoplus_{k=0}^n \bigwedge(n)_k$ 
and  $W(n) = \bigoplus_{k=-1}^{n-1} W(n)_k$, where 
\[
\bigwedge(n)_k:= \{ p(\xi_1,\ldots,\xi_n) \; | \; \deg p = k\}
\; \; \mbox{and} \; \;   
W(n)_{k}:=\{\sum_{l=1}^n p_i \frac{\partial}{\partial \xi_i} \; | \; \deg p_i = k+1\}.
\] 
In particular, $W(n)_0 \cong \gg \gl (n)$. Set 
$\bigwedge(n)_{\geq j} : = \bigoplus_{k \geq j} \bigwedge(n)_k$ and 
$W(n)_{\geq j} : = \bigoplus_{k \geq j} W(n)_k$. 
The Lie algebra $\gg_{\bar{0}}$ has the Levi subalgebra 
\[
\gl_{\bar{0}} = W(n)_0 \supset \gh_{\bar{0}}
\]
and nilradical
\[
\gn_{\bar{0}}^+= \bigoplus_{k \geq 1} W(n)_{2k}.
\]
The Weyl group $W_{\gl_{\bar{0}}}$ is isomorphic to the symmetric group $S_n$. 
Its action on the 
root lattice of $\gl_{\bar{0}}$ extends to actions on $\bigwedge(n)$ and 
$W(n)$ by Lie algebra automorphisms: for $\sigma \in S_n$, 
$\sigma (\xi_i) = x_{\sigma(i)}$, 
$\sigma (\partial/\partial \xi_i) = \partial/\partial \xi_{\sigma(i)}$.

For $1 \leq i \leq n$ denote
\[
\varepsilon_i \in \gh^*, \quad 
\varepsilon_i \left( \xi_l \frac{\partial}{\partial \xi_l} \right) = \delta_{il}.
\]
Given $I=\{ i_1, \ldots,i_k\} \subseteq [1,n]$ 
and $j \in [1,n]$ such that $j \notin I$, we set 
\[
\varepsilon_{I, j}:=\varepsilon_{i_1}+ \ldots + \varepsilon_{i_k} - \varepsilon_{j}.
\]
For $I=\{ i_1, \ldots,i_k\} \subsetneq [1,n]$, set
\[
\varepsilon_I:=\varepsilon_{i_1}+ \ldots + \varepsilon_{i_k}. 
\]
The root system of $\gg = W(n)$ is 
\[
\Delta = \{ \varepsilon_{I,j} \; | \; I \subseteq [1,n], j \in [1,n], j \notin I \}
\sqcup \{ \varepsilon_I \; | \; I \subsetneq [1,n] \}. 
\]
The corresponding root spaces are 
\[
\gg^{\varepsilon_{I, j} } = \Span \left\{ \xi_{i_1} \ldots \xi_{i_k} 
\frac{\partial}{\partial \xi_j} \right\}
\]
for $I=\{ i_1, \ldots,i_k\} \subseteq [1,n]$, $j \in ( [1,n] \backslash I )$  
and 
\[
\gg^{\varepsilon_{I} } = \Span \left\{ \xi_{i_1} \ldots \xi_{i_k} \xi_l 
\frac{\partial}{\partial \xi_l} \; \Big{|} \; 
l \in ([1,n] \backslash I) \right\}
\]
for $I=\{ i_1, \ldots,i_k\} \subsetneq [1,n]$.

Consider the subalgebra of $\gg = W(n)$
\[
\gs := \Span \left\{ \frac{\partial}{\partial \xi_i}, 
\xi_i \frac{\partial }{\partial \xi_j}, 
\xi_i \sum_{l=1}^n\xi_l\frac{\partial}{\partial \xi_l} 
\; \Big{|} \; 1\leq i \neq j \leq n 
\right\}
\supset \gh. 
\]
Its root system, considered as a subsystem of $\Delta$, is given by
\[
\Delta_{\gs} = \{ \varepsilon_i - \varepsilon_j \; | \; 1 \leq i \neq j \leq n  \}
\sqcup \{ \pm \epsilon_i \; | \; 1 \leq i \leq n \}.
\]  
We have the isomorphism
\begin{equation}
\label{isom}
\gs \cong \gs \gl (1 |n), \quad
\xi_i \frac{\partial }{\partial \xi_j} \mapsto 
E_{i+1, j+1} -\delta_{ij} E_{11}, \; 
\xi_i \sum_{l=1}^n\xi_l\frac{\partial}{\partial \xi_l} 
\mapsto E_{i+1, 1}, \; 
\frac{\partial}{\partial \xi_i} \mapsto E_{1, i+1},
\end{equation}
for $1 \leq i \neq j \leq n$.

Recall that the root system of $\gs \gl (1 |n)$ is 
\[
\{ \varepsilon_i - \varepsilon_j \; | \;  1 \leq i \neq j \leq n \} 
\sqcup \{ \pm(\delta_1 - \epsilon_i) \; | \; 
1 \leq i \leq n \}. 
\] 
(We interchanged the roles of $\varepsilon$ and $\delta$ from 
\S \ref{sl} in order to align that notation to the standard 
one for $W(n)$.) The root system of $\gs \gl (1 |n)$
is identified with $\Delta_{\gs}$ via the isomorphism
\eqref{isom} and $\delta_1 \mapsto 0$.

In order to describe the cominuscule parabolic sets of roots of $W(n)$,
we introduce the following sets
\begin{eqnarray*}
L_{W(n)}(n_0) &= &  \{ \varepsilon_{I} \; | \; I \subseteq [1,n_0] \} 
\sqcup \{ \varepsilon_{I, j} \; | \; 
I \subseteq [1,n_0], j \in [1,n_0], j \notin I\} 
\\
&& \sqcup \{ \varepsilon_{I \sqcup \{ i \}, j} \; | \; 
I \subseteq [1, n_0], i \neq j \in [n_0 +1, n] \},
\\
N^+_{W(n)}(n_0) &=& \{ \varepsilon_{I, j} \; | \; I \subseteq [1, n_0], 
j \in [n_0 +1, n] \},\\
N^-_{W(n)}(n_0) & =& \{ \varepsilon_I \; | \;  
I \subseteq [1,n], I \not\subseteq [1,n_0] \}
\\
&& \sqcup 
\{ \varepsilon_{I , j} \; | \;  
I \subseteq [1,n], I \not\subseteq [1,n_0], j \in [1, n_0], j \notin I \}
\\
&& \sqcup
\{ \varepsilon_{I , j} \; | \;  
I \subseteq [1,n], j \in [n_0 +1, n], j \notin I, |I \cap [n_0+1, n]| \geq 2 \},
\end{eqnarray*}
where $0 \leq n_0 < n$ and $|S|$ denotes the cardinality 
of a finite set $S$. The sets
\begin{equation}
\label{WP}
P_{W(n)}(n_0) = L_{W(n)}(n_0) \sqcup N_{W(n)}^+(n_0) \; \; 
\mbox{and} \; \;  
P^{-}_{W(n)}(n_0) = L_{W(n)}(n_0) \sqcup N^{-}_{W(n)}(n_0)
\end{equation}
are principal 
parabolic subsets of $\Delta$ with respect to the functionals
$\Lambda_{n_0}$ and $- \Lambda_{n_0}$, where
\begin{equation} \label{prin_w_n}
\Lambda_{n_0} (\varepsilon_i) = 0 \; \; 
\mbox{for} \; \; i \in [1,n_0], \; \; 
\Lambda_{n_0} (\varepsilon_i) = -1 \; \; 
\mbox{for} \; \; i \in [n_0+1, n].
\end{equation}
In particular,
$\Delta = N^+ (n_0) \sqcup L(n_0) \sqcup N^-(n_0)$
and the decompositions in \eqref{WP} are Levi decompositions.
Using Proposition \ref{algebras}, it is straightforward to verify 
that the sets $P_{W(n)}(n_0)$ are cominuscule, while the set 
$P^-_{W(n)} (n_0)$ is cominuscule if and only 
if $n_0 = n-1$. Denote by $\gl_{W(n)}(n_0)$
and $\gn^\pm_{W(n)}(n_0)$ the root 
subalgebras of $W(n)$ corresponding to the 
sets of roots $L_{W(n)}(n_0)$ and $N^\pm_{W(n)}(n_0)$,
and such that $\gl_{W(n)} (n_0) \supset \gh$ and 
$\gn^\pm(n_0) \cap \gh = 0$, recall \eqref{root-sub}. 
For $1 \leq k \leq l \leq n$ denote the subalgebra
\[
\gg \gl [k,l] = \Span \left\{ 
\xi_i \frac{\partial }{\partial \xi_j} \; \Big{|} \; 
k \leq i, j \leq l \right\}
\]
of $W(n)$. Clearly $\gg \gl [k,l] \cong \gg \gl (l-k)$.

\begin{lemma} 
\label{Wn-alg} For all $0 \leq n_0 \leq n-1$  
we have the isomorphism of Lie superalgebras
\[
\gl_{W(n)}(n_0) \cong 
W(\xi_1, \ldots,\xi_{n_0}) \ltimes 
\left( \bigwedge(\xi_1, \ldots ,\xi_{n_0}) \otimes \gg \gl [n_0+1, n] \right),
\]
where the second term represents the Lie superalgebra
which is the tensor product of a supercommutative 
algebra and a Lie superalgebra, and  
$W(\xi_1, \ldots,\xi_{n_0})$ acts on 
$\bigwedge(\xi_1, \ldots ,\xi_{n_0})$ by derivations. Moreover, 
we have the isomorphisms of $\gl_{W(n)}(n_0)$-modules
\begin{eqnarray*}
\gn^+_{W(n)}(n_0) & \cong & \bigwedge (\xi_1, \ldots ,\xi_{n_0}) \otimes (V^{n-n_0})^* 
\; \; \mbox{and} \\
\gn^-_{W(n)}(n_0) & \cong & W(\xi_1, \ldots,\xi_{n_0})   \otimes 
\bigwedge(\xi_{n_0+1}, \ldots, \xi_n)_{\geq 1}  \\
&&
\oplus \bigwedge(\xi_1, \ldots, \xi_{n_0}) 
\otimes W(\xi_{n_0+1}, \ldots, \xi_n)_{\geq 1}.
\end{eqnarray*}
For the module structure of the right hand sides we use 
the adjoint actions of the chain of Lie subalgebras
$\gg \gl [n_0+1, n] \subset W(\xi_{n_0+1}, \ldots, \xi_n)_{\geq 1} 
\subset W(\xi_{n_0+1}, \ldots, \xi_n)$ and the (left)
multiplication action of the supercommutative algebra 
$\bigwedge(\xi_1, \ldots ,\xi_{n_0})$ on $W(\xi_1, \ldots,\xi_{n_0})$
and itself. The symbol $(V^{n-n_0})^*$ denotes the dual of the 
vector representation of $\gg \gl [n_0+1, n]$.
\end{lemma}

The proof of of Lemma \ref{Wn-alg} amounts to a direct computation 
of the root algebras $\gl_{W(n)}(n_0)$ and $\gn^\pm_{W(n)}(n_0)$ 
and is left to the reader. The next result classifies and 
describes the cominuscule parabolic sets of roots of $W(n)$.

\begin{theorem} 
\label{prop_wn}
There are $n+1$ orbits of cominuscule parabolic subsets of the root system 
of $W(n)$ under the action of the Weyl group $W_{\gl_{\bar{0}}}\cong S_n$.
The parabolic sets $P = P_{W(n)}(n_0)$, $0 \leq n_0 \leq n-1$, and 
$P = P_{W(n)}^{-}(n-1)$ provide representatives of those orbits.
They have unique Levi decompositions given by
\begin{eqnarray}
\label{WnLevi1} 
P_{W(n)}(n_0) &=& L_{W(n)}(n_0) \sqcup N_{W(n)}(n_0), \;
0 \leq n_0 \leq n-1 \; \; \mbox{and} 
\\
P_{W(n)}^-(n-1) &=& L_{W(n)}(n-1) \sqcup N_{W(n)}^-(n-1).
\label{WnLevi2}
\end{eqnarray}
The Levi components $\gl_{W(n)}(n_0)$ of the corresponding 
parabolic subalgebras of $W(n)$ and their nilradicals $\gn^\pm_{W(n)}(n_0)$
considered as $\gl_{W(n)}(n_0)$-modules are given by Lemma \ref{Wn-alg}
\end{theorem}
\noindent
{\bf Proof.} Assume that $P$
is a cominuscule parabolic subset of $\Delta$,
and that $P = L \sqcup N^+$ is a 
Levi decomposition of $P$.
First we show that $\Delta_{\gs} \not\subseteq P$. Indeed,
if $\Delta_{\gs} \subseteq P$, then $\pm \varepsilon_i \in L$ 
for all $i \in [1,n]$. This implies that $\Delta = P$, which 
contradicts to the condition that $P$ is a proper subset of $\Delta$.

Therefore by Lemma \ref{restrict}, $P \cap \Delta_{\gs}$ is 
a cominuscule parabolic subset of $\Delta_{\gs}$. Observe that 
$\gs_{\bar{0}}= \gl_{\bar{0}}$, thus 
$W_{\gs_{\bar{0}}} \cong W_{\gl_{\bar{0}}} \cong S_n$.
Using Theorem \ref{prop_slmn} and 
the isomorphism \eqref{isom} we obtain that
there exist integers 
$0 \leq m_0 \leq 1$ and $0 \leq n_0 \leq n$, 
$(m_0, n_0) \neq (0,0)$, $(1,n)$, such that $P$ is conjugated 
under the action of $W_{\gl_{\bar{0}}}$ to a cominuscule 
parabolic subset of $\Delta$ such that
\begin{equation}
\label{LNW}
L \cap  \Delta_{\gs} =  L_{\gs } (m_0, n_0)
\; \; \mbox{and} \; \; 
N^+ \cap  \Delta_{\gs} = N_{\gs }^+ (m_0, n_0),
\end{equation}
where
\begin{eqnarray}
\label{LW}
L_{\gs } (m_0, n_0) & = & 
\{ \varepsilon_i - \varepsilon_j \; | \; 
i, j \leq n_0 \; \mbox{or} \; i, j > n_0 \}, 
\\
&& \sqcup \{ \pm \varepsilon_i 
\; | \; i \leq n_0 \; \mbox{if} \; m_0 = 1, 
i> n_0 \; \mbox{otherwise} \} 
\nonumber
\\
\label{NW0}
N_{\gs }^+ (0, n_0)
& = & 
\{ \varepsilon_i - \varepsilon_j \; | \; i \leq n_0 < j \} 
\sqcup \{ \varepsilon_i \; | \; i \leq n_0 \},
\\
N_{\gs }^+ (1, n_0)
& = & 
\{ \varepsilon_i - \varepsilon_j \; | \; 
i \leq n_0 < j \} \sqcup \{ - \varepsilon_j \; | \; 
j >  n_0 \}.
\label{NW1}
\end{eqnarray}
We conjugate $P$ by an element of $W_{\gl_{\bar{0}}}$ so that 
\eqref{LNW} holds.

We first consider the case $m_0 = 0$. If $n_0 > 1$ then 
from \eqref{LNW} and \eqref{NW0} we obtain 
$\varepsilon_1, \varepsilon_2 \in N^+$.
Since $\varepsilon_1 + \varepsilon_2 \in \Delta$, 
Proposition \ref{algebras}
leads to a contradiction. Thus $n_0 =1$. 
Lemma \ref{sums} implies that 
\begin{eqnarray}
L &\supseteq& 
\{ \varepsilon_{I, j} \; | \; I \subseteq [2,n], j \in [2, n], 
j \notin I \} \sqcup \{ \varepsilon_I \; | \; I \subseteq [2,n] \},
\label{Lin1}
\\
N^+ &\supseteq&
\{ \varepsilon_{ \{ 1 \} \sqcup I, j} \; | \; I \subseteq [2,n], j \in [2, n], 
j \notin I \} \sqcup \{ \varepsilon_{ \{1 \} \sqcup I} 
\; | \; I \subseteq [2,n] \},
\label{Nin1}
\end{eqnarray}
where $w_0$ is the longest element of $W_{\gl_{\bar{0}}} \cong S_n$.
If any of these inclusions are strict, then $P$ contains an 
element of the form $\varepsilon_{I,1}$, where $I \subseteq [2,n]$, 
$1 \notin I$. Since $\pm\epsilon_j \in P$ for all $j \in [2,n]$, this 
would imply that $- \epsilon_1 \in P$. Therefore $P = \Delta$, 
because $\epsilon_1 \in P$. This is a contradiction. Thus 
both inclusion \eqref{Lin1} and \eqref{Nin1} are 
equalities, which implies that   
$w_0 (L)=L_{W(n)}(n-1)$ and $w_0 (N^+) =N_{W(n)}^{-}(n-1)$, where 
$w_0$ is the longest element of $W_{\gl_{\bar{0}}} \cong S_n$.

Now let $m_0=1$.  We will show that $L = L_{W(n)}(n_0)$ 
and $N^+ = N^+_{W(n)}(n_0)$. First we prove that
\begin{equation}
\label{WWin}
L \supseteq L_{W(n)}(n_0) \; \; 
\mbox{and} 
\; \;  N^+ \supseteq N^+_{W(n)}(n_0).
\end{equation}
Since for all $i \in [1, n_0]$, $\pm \varepsilon_i \in L$, 
Lemma \ref{sums} (iii) implies that for all $I \subseteq [1, n_0]$, 
$\varepsilon_I \in L$, and for all $I \subseteq [1, n_0]$ and $j \in [1,n_0]$,
$j \notin I$, $\varepsilon_{I, j} = \varepsilon_I - \epsilon_j \in L$.
Similarly for all $I \subseteq [1, n_0]$ and $i \neq j \in [n_0 +1, n]$,
$\varepsilon_{I \sqcup \{ i\}, j} = \varepsilon_I + 
(\varepsilon_i - \epsilon_j) \in L$. This proves the first 
inclusion in \eqref{WWin}.

Since $\varepsilon_I \in L$ for $I \subseteq [1,n_0]$ and 
by \eqref{NW1} $- \epsilon_j \in N^+$ for $j \in [n_0+1, n]$, 
Lemma \ref{sums} (ii) implies that $\varepsilon_{I, j} \in N^+$, 
under the same conditions on $I$ and $j$. This proves the second
inclusion in \eqref{WWin}.

To prove that the inclusions in \eqref{WWin} are equalities, 
we need to show that
\begin{equation}
\label{inter}
P \cap N^-_{W(n)} (n_0) = \emptyset.
\end{equation}

First we show that 
\begin{equation}
\label{first}
\varepsilon_I \notin P, \; \; 
\mbox{for} \; \;  
I \subseteq [1,n], I \not\subseteq [1,n_0]. 
\end{equation}
Assume 
the opposite. Then $-\varepsilon_j \in N^+$, 
$\forall j \in I \cap [n_0+1, n]$
and Lemma \ref{sums} (ii), (iv) imply that
$\varepsilon_{I \cap [1,n_0]} = \varepsilon_I +
\sum_{j \in I \cap [n_0+1, n]} (-\varepsilon_j) \in N^+$. 
Since we already showed that 
$\varepsilon_{I \cap [1,n_0]} \in L$, 
this contradicts with $L \cap N^+ = \emptyset$.
  
Next we prove that $\varepsilon_{I,j} \notin P$ for 
$I \subseteq [1,n]$, $I \not\subseteq [1,n_0]$, 
$j \in [1, n_0]$. If this is not the case, then 
$\varepsilon_j \in L$, 
$\forall j \in [1,n_0]$
and \eqref{P1}
imply that 
$\varepsilon_I = \varepsilon_{I,j} + \varepsilon_j \in P$,
which contradicts with \eqref{first}.

Finally, we prove that $\varepsilon_{I, j} \notin P$ 
for $I \subseteq[1,n]$, $j \in [n_0+1,n]$ such that 
$j \notin I$ and $|I \cap [n_0+1, n]| \geq 2$. 
Assume the opposite and choose 
$i \in I \cap[n_0+1, n]$. Then $\varepsilon_j -\varepsilon_i \in L$, 
$\forall j \in [n_0 +1, n]$ and \eqref{P2} imply that 
$\varepsilon_{ I \backslash \{ i \} } = 
\varepsilon_{I, j} + (\varepsilon_j - \varepsilon_i)
\in P$, which again contradicts with \eqref{first}
since $I \backslash \{ i \} \not\subseteq [1,n_0]$.

This proves that each cominuscule parabolic set of roots for $W(n)$ 
is conjugated under the action of the Weyl group 
$W_{\gl_{\bar{0}}}$ to one of the sets 
\begin{equation}
\label{sets}
P = P_{W(n)}(n_0), \; 0 \leq n_0 \leq n-1, \; \; 
\mbox{and} \; \; 
P = P_{W(n)}^{-}(n-1),
\end{equation}
and that those parabolic subsets are principal 
and have unique Levi decompositions 
given by \eqref{WnLevi1}--\eqref{WnLevi2}. 
It remains to show that none of those parabolic sets 
are in the same $W_{\gl_{\bar{0}}}$-orbit. 
Since $\gl_{\bar{0}} = \gs_{\bar{0}}$, the set of 
roots $\Delta_{\gs}$ is stable under the action 
of $W_{\gl_{\bar{0}}}$. If two parabolic sets 
of roots of $W(n)$, $P$ and $P'$ are conjugated under 
$W_{\gl_{\bar{0}}}$, then 
$P' \cap \Delta_{\gs} \in W_{\gl_{\bar{0}}} 
(P \cap \Delta_{\gs})$. Since
\begin{eqnarray*}
P_{W(n)}(n_0) \cap \Delta_{\gs} 
&=& L_{\gs} (1, n_0) \sqcup N_{\gs}^+ (1, n_0) , \; 
1 \leq n_0 \leq n-1, 
\\
P_{W(n)}^-(n-1) \cap \Delta_{\gs} 
&=& w_0( L_{\gs} (1, n-1) \sqcup N_{\gs}^+ (0, n-1)),
\end{eqnarray*}
where $w_0$ is the longest element of $W_{\gl_{\bar{0}}} \cong S_n$,
Theorem \ref{prop_slmn} implies that none of the parabolic sets 
\eqref{sets} are in the same $W_{\gl_{\bar{0}}}$-orbit.
\hfill $\square$
\medskip

\begin{remark}
Following the proof of Theorem \ref{prop_wn} we may define alternatively the sets $L_{W(n)}(n_0)$,  $N_{W(n)}(n_0)$,  and $N^{-}_{W(n)}(n_0)$ as the unique sets $L, N^{+}, N^{-}$ that satisfy the properties of Lemma \ref{roots} and such that $\pm \varepsilon_i \in L$, $1 \leq i \leq n_0$, $\varepsilon_j \in N^{-}$, $- \varepsilon_j \in N^{+}$, $n_0 < j \leq n$. 
\end{remark}

The proof of Theorem \ref{prop_wn} implies the following result.
\begin{corollary} \label{cor_wn}
If $n_0 > 1$, then the cominuscule parabolic  subset of roots 
$P_{\gs \gl (n_0)} (n_0)$ of $W(n)_0$ has a unique cominuscule 
parabolic $W(n)$-extension: $P_{W(n)}(n_0)$. The cominuscule parabolic 
subset  $P_{\gs \gl (n_0)} (1)$ has two distinct 
extensions: $P_{W(n)}(1)$ and $w_0 P^{-}_{W(n)}(n-1)$. The cominuscule parabolic 
subsets of roots $P_{\gs \gl (1|n)}(1|n_0)$ 
of $\gs \cong \gs \gl (1|n)$ have no cominuscule parabolic $W(n)$-extensions 
for $n_0>1$. In all other cases, $P_{\gs \gl (1|n)}(m_0|n_0)$ has a unique 
cominuscule parabolic 
$W(n)$-extension $\widehat{P}_{\gs \gl (1|n)}(m_0|n_0)$: 
$\widehat{P}_{\gs \gl (1|n)}(0|1) = P^{-}_{W(n)}(n-1)$ and  
$\widehat{P}_{\gs \gl (1|n)}(1|n_0) = P_{W(n)}(n_0)$ for $0 \leq n_0 \leq n-1$. 
\end{corollary}

\begin{remark}
There are two important particular cases in Theorem \ref{prop_wn}. The first such case is $n_0 = 0$ when we have 
$\gl = W(n)_0$, $\gn^+ = W(n)_{-1}$ and 
$\gn^{-} = W(n)_{\geq 1}$. The other case is $n_0 = n-1$ which corresponds 
to $\gl \cong W(n-1) \ltimes \left( \gg \gl (1) \otimes \bigwedge (n-1) \right)$, and $\gn^{+}$ and $\gn^{-}$ 
are isomorphic to  $\bigwedge (n-1) \otimes V^*$ and $W (n-1) \otimes V$, respectively, where $V$ is the standard one dimensional
$\gg \gl (1)$-module.
\end{remark}

Theorem \ref{prop_wn} also implies:

\begin{remark}
All cominuscule parabolic sets of roots of $W(n)$ 
are principal parabolic subsets with functionals 
defined in (\ref{prin_w_n}).
\end{remark}
\subsection{$\gg = S(n)$ and $\gg = S'(n)$.}
\label{S_n}
There are two associative superalgebras of differential forms defined over $\bigwedge(n)$, namely $\Omega(n)$ and $\Theta(n)$. 
The superalgebra $\Omega(n)$ has generators $d \xi_1, \ldots,d \xi_n$ and 
defining relations $d \xi_i \circ d \xi_j = d \xi_j  \circ d \xi_i$, $\deg d\xi_i = \bar{0}$, while the superalgebra 
$\Theta(n)$ has generators  $\theta \xi_1, \ldots,\theta \xi_n$ and relations 
$\theta  \xi_i \wedge \theta \xi_j = - \theta \xi_j  \wedge \theta \xi_i, \deg \theta \xi_i = \bar{1}$. Note that the 
differentials $d$ and $\theta$ are  derivations of degree $\bar{1}$ and  $\bar{0}$ respectively. Let 
$\mu_n:= \theta \xi_1 \wedge \ldots \wedge \theta \xi_n$ be the standard volume form in $\Omega(n)$ and $\mu_n':=(1 + \xi_1 \ldots \xi_n) \mu_n$.

Every derivation $D$ of $W(n)$ and every automorphism $\Phi$ of $\bigwedge(n)$ extend uniquely to 
a derivation $\widetilde{D}^{d}$ and an automorphism $\widetilde{\Phi}^{d}$ 
(respectively, $\widetilde{D}^{\theta}$ and $\widetilde{\Phi}^{\theta}$) of $\Omega(n)$ 
(resp., $\Theta(n)$) so that $[\widetilde{D}^{d}, d]= 0$ and $[\widetilde{\Phi}^{d}, d]= 0$  
(resp., $\widetilde{D}^{\theta} \theta f - \theta  \widetilde{D}^{\theta}f =0$ and  
$\widetilde{\Phi}^{\theta} \theta f - \theta  \widetilde{\Phi}^{\theta}f =0$  for every 
$f \in \bigwedge(n)$). We denote by $S(n)$ 
the Lie superalgebra $\{ D \in W(n)\; | \; \widetilde{D}^\theta(\mu_n) = 0\}$ 
and by $S'(n)$ the Lie superalgebra  $\{ D \in W(n)\; | \; \widetilde{D}^\theta(\mu_n') = 0\} $. 
Since $S'(2k+1) \cong S(2k+1)$ we consider $S'(n)$ only for even numbers 
$n$. In explicit terms we have:
\begin{eqnarray*}
S(n) & = & \Span \left\{ \frac{\partial f}{\partial \xi_i} \frac{\partial}{\partial \xi_j} 
+  \frac{\partial f}{\partial \xi_j} \frac{\partial}{\partial \xi_j}\; | \; f \in \bigwedge(n), 
1 \leq i,j \leq n\right\},\\
S'(n) &  = & \Span \left\{ (1 - \xi_1 \ldots \xi_n)\left(\frac{\partial f}{\partial \xi_i} 
\frac{\partial}{\partial \xi_j} +  \frac{\partial f}{\partial \xi_j} \frac{\partial}{\partial \xi_j}\right)\; 
| \; f \in \bigwedge(n), i,j = 1, \ldots,n\right\}.
\end{eqnarray*}

In particular, $S(n)_0 = S'(n)_0 \cong \gs \gl (n)$. Set $S(n)_{\geq j}:= \bigoplus_{i \geq j} S(n)_i$. 
Note that $S'(n)$ is not a graded Lie subalgebra of $W(n)$, but it has a filtration induced by the 
filtration $\{ W(n)_{\geq j}\}_j$ of $W(n)$. The corresponding graded superalgebra is isomorphic to $S(n)$.

We fix the Cartan subalgebra of $S(n)$ to be $\gh_{S(n)} = \gh_{W(n)} \cap S(n)$ 
and set $\gh_{S'(n)}  = \gh_{S(n)} $ for even $n$. 
The root systems of $S(n)$ and $S'(n)$ coincide and can be described as follows. 
Denote by $\iota_S : S(n) \to W(n)$ the natural inclusion. Let $ \overline{\Delta}_{S(n)}$ be  
the set obtained from ${\Delta}_{W(n)}$ by removing 
the $n$ roots $\varepsilon_{[1,n] \backslash \{ i\} }$, $1 \leq i \leq n$:
\[
\overline{\Delta}_{S(n)} =  \{ {\varepsilon}_{I,j} \; | \; I \subseteq [1,n], j \in [1,n], j \notin I \}
\sqcup \{ {\varepsilon}_I \; | \; | I |  \leq n-2  \}. 
\]
The kernel of the restriction map $ {\iota_S^* |}_{ \gh_{W(n)}^*} : \gh_{W(n)}^* \to \gh_{S(n)}^*$  
equals $\C (\varepsilon_1 + \ldots + \varepsilon_n)$. Throughout this subsection,
we will identify the root system $\Delta_{S(n)}$ with $\overline{\Delta}_{S(n)}$
via $\iota_S^*$. In particular, by abuse of notation, for a root 
$\alpha \in \overline{\Delta}_{S(n)}$, we will write $\alpha$ instead of $\iota_S^*(\alpha)$.
All sums of roots of $\Delta_{S(n)}$ will be computed in $\gh_{W(n)}^*$
via this identification.

The root spaces of $S(n)$ are described as follows:
\[
S(n)^{{\varepsilon}_{I, j} } = \Span \left\{ \xi_{i_1} \ldots \xi_{i_k} 
\frac{\partial}{\partial \xi_j} \right\}
\]
for $I=\{ i_1, \ldots,i_k\} \subseteq [1,n]$, $j \notin I$  
and 
\[
S(n)^{{\varepsilon}_{I} } = \Span \left\{ \xi_{i_1} \ldots \xi_{i_k}\left(  \xi_l 
\frac{\partial}{\partial \xi_l} -  \xi_m 
\frac{\partial}{\partial \xi_m} \right) \; \Big{|} \; 
l,m \in ([1,n] \backslash I) \right\}
\]
for $I=\{ i_1, \ldots,i_k\} \subsetneq [1,n]$.

On the other hand, the root spaces of $S'(n)$ coincide with those of $S(n)$ except 
for $S'(n)^{-{\varepsilon }_j}$, $j=1,\ldots ,n$. For the latter root spaces we have
\[
S'(n)^{-{\varepsilon}_{j} } = \Span \left\{ (1 - \xi_1 \ldots \xi_n)
\frac{\partial}{\partial \xi_j} \right\}.
\]

The following lemma can be proved using the explicit description of $\gg^{\alpha}$ with the same reasoning 
as in the proof of Proposition \ref{appendix2} for $\gg = W(n)$.
\begin{lemma} \label{lemma_sn}
Let $\gg = S(n)$ or $\gg = S'(n)$, and let $\alpha, \beta \in {\Delta}_{S(n)}$ be such 
that ${\alpha} + {\beta} \neq 0$. Then

(i) $\left[ \gg^{\alpha},  \gg^{\beta}\right] \neq 0$ if and only 
if $\alpha + \beta \in {\Delta}_{W(n)}$ for $\gg = S(n)$;

(ii) $\left[ \gg^{\alpha},  \gg^{\beta} \right] \neq 0$ if and only 
if $\alpha + \beta \in {\Delta}_{W(n)}$ or $\alpha = - {\varepsilon}_i$, 
$\beta = - {\varepsilon}_j$, some $i \neq j$, for $\gg = S'(n)$.
\end{lemma}

Recall that we identify ${\Delta}_{S(n)}$ and $\overline{\Delta}_{S(n)}$ and the sum $\alpha + \beta$ 
in (i) and (ii) is taken in $\gh_{W(n)}^*$.

We define $L_{S(n)} (n_0)$, $N^{+}_{S(n)} (n_0)$, $N^{-}_{S(n)} (n_0)$, $P_{S(n)} (n_0)$, 
and $P^{-}_{S(n)} (n_0)$ with the same formulas that we used for the corresponding sets in \S \ref{W_n}. 

Contrary to $W(n)$, $S(n)$ does not have a subalgebra whose root system is 
$\{ \pm \varepsilon_i,  \varepsilon_i -  \varepsilon_j  \; | \; i\neq j \}$. However, 
we may define a monomorphism $\iota_W : W(n) \to S(n+1)$ by 
\[ 
\xi_{i_1} \ldots \xi_{i_k} \frac{\partial}{\partial \xi_j} \mapsto 
\xi_{i_1} \ldots \xi_{i_k}  
\frac{\partial}{\partial \xi_j}
\; \;  \mbox{and} \; \; 
\xi_{i_1} \ldots \xi_{i_k} \xi_j \frac{\partial}{\partial \xi_j} 
\mapsto\xi_{i_1} \ldots \xi_{i_k}  \left(\xi_j \frac{\partial}{\partial \xi_j} -  
\xi_{n+1}\frac{\partial}{\partial \xi_{n+1}}\right),
\]  
for $j \notin \{ i_1, \ldots ,i_k\}$.  

The Lie algebra $\gl_{\bar{0}}:=S(n)_0 =  S'(n)_0$ is a Levi subalgebra 
of $S(n)_{\bar{0}}$ and $S'(n)_{\bar{0}}$ containg  
$\gh_{S(n)}$ and $\gh_{S'(n)}$. The Weyl group $W_{\gl_{\bar{0}}}$ is 
isomorphic to the symmetric group $S_n$. It acts on $S(n)$ and $S'(n)$ 
by Lie algebra automorphisms and on the corresponding root systems, 
as follows. The automorphisms $\sigma \in S_n$ from \S \ref{W_n}
leave invariant the volume forms $\mu_n$ and $\mu_n'$. Thus
each $\sigma \in S_n$ preserves $S(n)$ and $S'(n)$.
The corresponding action of $S_n$ on the root systems
of $S(n)$ and $S'(n)$ is simply the restriction of the 
action of $S_n$ from $\Delta_{W(n)}$ to $\overline{\Delta}_{S(n)}$.   
The following theorem classifies the cominuscule 
parabolic subsets of roots of $S(n)$ and shows 
that each of them can be obtained by   
restricting a cominuscule parabolic subset 
of roots of $W(n)$. 

\begin{theorem} 
\label{prop_sn}
(i) There are $n+1$ orbits of cominuscule parabolic subsets of the root system 
of $S(n)$ under the action of the Weyl group $W_{\gl_{\bar{0}}}\cong S_n$. The parabolic subsets 
$P_{S(n)} (n_0)$, $0 \leq n_0 \leq n-1$, and $P = P^{-}_{S(n)} (n-1)$ provide representatives of these orbits. 
Each cominuscule parabolic subalgebra of $S(n)$ has a unique Levi decomposition.
The Levi subalgebra $\gl$ and the niradicals $\gn^+$ and $\gn^-$ of the 
above cominuscule parabolic subalgebras can be obtained by intersecting the corresponding 
subalgebras of $W(n)$ described in Theorem \ref{prop_wn} with $S(n)$. 
 
(ii) The Lie superalgebra $S'(n)$ has no cominuscule parabolic subsets.
\end{theorem}
\noindent
{\bf Proof.} We start with $\gg = S(n)$. Let $P$ be a cominuscule 
parabolic subset of $\Delta$ and $P= L \sqcup N^+$ 
be a Levi decomposition of $P$ for which the coresponding niradical
$\gn^+$ of $\gp_P$ is abelian, recall Definition \ref{LN-alg}.
Let $N^-$ be as in Lemma \ref{sums}.

For simplicity of the notation in this proof we set $\Delta = \Delta_{S(n)}$ and 
$\Delta_{\gs \gl(n)} = \Delta_{S(n)_0}$.  Assume first that $N^+ \cap \Delta_{\gs \gl (n)} \neq  \emptyset$. 
Then by Proposition \ref{algebras}, there exits $n_0$, $1\leq n_0 \leq n-1$, 
such that after conjugating $P$ by an element $\sigma \in S_n$ we have
$L\cap \Delta_{\gs \gl (n)}  = L_{\gs \gl (n)}(n_0)$ and 
$N^+ \cap \Delta_{\gs \gl (n)}  = N_{\gs \gl (n)}^+(n_0)$. We proceed with 
a case-by-case verification using Lemma \ref{lemma_sn}.

{\it Case 1: $\varepsilon_{1} \in N^+$}. 
Then $\varepsilon_{j} = \varepsilon_1 + (\varepsilon_j- \varepsilon_1) \in N^+$ 
for every $j \leq n_0$. In particular, $P$ is not cominuscule if $n_0 > 1$
since $\varepsilon_1 + \varepsilon_j \in \Delta_{W(n)}$, cf. Lemma \ref{lemma_sn}.
If $n_0 = 1$, 
then $\varepsilon_j \in L$ for every $j >1$. Indeed, if $\varepsilon_j \notin L$, then 
either $\varepsilon_j \in N^+$ or $-\varepsilon_j \in N^+$, which together with $\varepsilon_1 \in N^+$ leads to a contradiction since 
$\varepsilon_1 \pm \varepsilon_j \in \Delta_{W(n)}$. 
From here it is not difficult to verify that $P = w_0 P^{-}_{S(n)}(n-1)$.

{\it Case 2: $\varepsilon_{1} \in N^-$}. For $j \leq n_0$, 
$(\varepsilon_j- \varepsilon_1) \in L$, and for $j > n_0$, 
$(\varepsilon_j- \varepsilon_1) \in N^-$. Thus for all $j$,
$\varepsilon_j  = \varepsilon_1 + (\varepsilon_j- \varepsilon_1) \in N^-$
and $-\varepsilon_j \in N^+$. It is easy to conclude from here 
that $P = P_{S(n)}(0)$, which is a contradiction to 
$N^+ \cap \Delta_{\gs \gl (n)} \neq \emptyset$. 

{\it Case 3: $\varepsilon_1 \in L$}. Now we have 
$\varepsilon_j = \varepsilon_1 + (\varepsilon_j- \varepsilon_1) \in L$ 
for every $j \leq n_0$. If $\varepsilon_n \in L$, then we  
$\varepsilon_i = \varepsilon_n + (\varepsilon_i - \varepsilon_n) \in L$
for all $i > n_0$. Thus $\varepsilon_n \in L$ implies
$\varepsilon_j \in L$ for all $j$ 
and hence $L = P = \Delta$, which is a contradiction.
If $\varepsilon_n \in N^+$, then  $\varepsilon_1 =  \varepsilon_n  + 
(\varepsilon_1 - \varepsilon_n) \in N^+$, 
which is again a contradiction. 
It remains to consider the case when $\varepsilon_n \in N^-$. 
Then $ \varepsilon_j \in N^-$ for $j > n_0$. From here 
one easily verifies that $P = P_{S(n)}(n_0)$. 

Now we assume $N^+ \cap \Delta_{\gs \gl (n)} =  \emptyset$, 
and hence $\Delta_{\gs \gl (n)} \subset L$. 
Then either all $\varepsilon_j$ are in $N^+$ or all $\varepsilon_j$ are in $N^-$. 
(Otherwise there will exist two indices $i \neq j$ such that 
$\varepsilon_i, - \varepsilon_j \in N^+$, which contradicts to 
the assumption that $P$ is cominuscule 
since $\varepsilon_i - \varepsilon_j \in \Delta_{W(n)}$, 
cf. Lemma \ref{lemma_sn}.) Using once again the assumption 
that $P$ is cominuscule and Lemma \ref{lemma_sn} 
we rule out the case $\varepsilon_j \in N^+$, $\forall j$ since 
$\varepsilon_i + \varepsilon_j \in \Delta_{W(n)}$, 
$\forall i \neq j$. Therefore $\varepsilon_j \in N^-$
for all $j$, which leads to $P = P_{S(n)}(0)$.  

The isomorphisms for the Levi subalgebras and nilradicals of the cominuscule 
parabolic subalgebras obtained in this way are analogous to the $W(n)$ case.
The case $\gg = S'(2l)$ follows from $\gg = S(2l)$ and is left to the reader. 
\hfill $\square$

\begin{remark}
The explicit isomorphisms for $\gl$, $\gn^+$, and $\gn^-$ are analogous to the 
$W(n)$ case but are rather lengthy and will be omitted. In the particular 
case $n_0 = 0$, we have
$\gl = S(n)_0$, $\gn^+ = S(n)_{-1}$ and $\gn^{-} = S(n)_{\geq 1}$, while for $n_0 = n-1$, 
$\gl = \iota_W(W(\xi_1, \ldots ,\xi_{n-1}))$, $\gn^+ \cong \bigwedge 
(\xi_1, \ldots ,\xi_{n-1}) \otimes V^*$ and 
$\gn^{-} \cong S (\xi_1, \ldots ,\xi_{n-1}) \otimes V$, where $V$ denotes 
the standard $\gg \gl [n-1,n]$-representation. 
\end{remark}

\begin{remark}
One can prove Theorem \ref{prop_wn} using the same reasoning as 
in the proof of Theorem \ref{prop_sn}. 
The advantage of the present proof of Theorem \ref{prop_wn} is that it 
provides a valuable connection between 
the cominuscule parabolics of $\gs \gl (1|n)$ and $W(n)$ 
(cf. Corollary \ref{cor_wn}) which is not obvious otherwise.
\end{remark}

\begin{remark}
All cominuscule parabolic sets of roots of $S(n)$ 
are principal parabolic subsets with functionals 
defined in (\ref{prin_w_n}).
\end{remark}
\subsection{$\gg = H(n)$.}
\label{H_n}
The Hamiltonian finite dimensional Lie superalgebras are defined by 
$\widetilde{H}(n):= \{ D \in W(n)\; | \; \widetilde{D}^d\omega_n = 0\}$ 
and $H(n):=[\widetilde{H}(n), \widetilde{H}(n)]$. In explicit form: 
\begin{eqnarray*}
\widetilde{H}(n)& = & \Span \left\{
D_f:=\sum_i \frac{\partial f}{\partial \xi_i} 
\frac{\partial}{\partial \xi_i} \; | \; f \in
\bigwedge(n), f(0)=0, \; 1 \leq i,j \leq n \right\},\\
\widetilde{H}(n) & = & H(n) \oplus \C D_{\xi_1 \ldots \xi_n}.
\end{eqnarray*}
If we consider $\bigwedge(n)$ as a Poisson superalgebra, 
then the map ${\mathcal D} : \bigwedge(n) \to \widetilde{H}(n)$, 
$f \mapsto D_f$, is a surjective homomorphism of Lie superalgebras 
with $\ker {\mathcal D}  = \C$. In particular, $[D_f,D_g] = D_{\{f,g\}}$ 
where $\{f,g\}:=(-1)^{\deg f}\sum_{i=1}^n \frac{\partial f}{\partial \xi_i} 
\frac{\partial g}{\partial \xi_j}$.

The Lie superalgebra $H(n)$ is a graded subalgebra of $W(n)$. Set 
$H(n)_k = H(n) \cap W(n)_k$, $-1 \leq k \leq n-1$. 
We have $H(n)_0 \cong \mathfrak{so}(n)$. Every Cartan subalgebra $\gh$ of $H(n)$ has a nilpotent part. An explicit description of such subalgebras can be found in Appendix A of \cite{GP2}. We fix such $\gh$ for which $\gh \cap H(n)_0$ equals $ \Span \left\{D_{\xi_i \xi_{i+l}} \; | \; i=1, \ldots , [n/2] \right\}$.
The root system of $H(n)$ is given by
$\Delta = \{ \varepsilon_{I} - \varepsilon_{J}  \; | \; I, J \subset
[1, [ n/2 ] ], I \cap J = \emptyset \}$, where 
$\varepsilon_I = \sum_{i \in I} \varepsilon_i$ and the 
arithmetic of $\varepsilon_i$'s is the same as in the case of $W(n)$. 
All roots vanish on $\gh'_{\bar{0}}$.
Denote by $\Delta_{\gs \go (n)}$ the subset of roots corresponding to 
$H(n)_0$.

The Lie algebra $H(n)_0$ is a Levi subalgebra of $\gg_{\bar{0}}$.
As in the previous two subsections, every element $w$ of the Weyl group of 
$H(n)_0$ can be extended to a Weyl automorphism $\sigma_w$ of $H(n)$. 
This induces an action of $W_{H(n)}$ on $\Delta$ and 
on the set of parabolic subsets of $\Delta$. 
We define
\begin{eqnarray*}
L_{H(n)} &:= &  \{ \varepsilon_I -  \varepsilon_J \; | \; 
I, J \subset [1, [ n/2 ] ], I \cap J = \emptyset, 
1\notin I, 1 \notin J\},\\
N_{H(n)}^{+} &:=&  \{ \varepsilon_I -  \varepsilon_J \; | \; 
I, J \subset [1, [ n/2 ] ], I \cap J = \emptyset, 
1 \in I\},\\
\end{eqnarray*}
and $P_{H(n)} :=L_{H(n)} \sqcup N_{H(n)}^{+}$.
One easily checks that $\Delta = (- N_{H(n)} ^+) \sqcup L_{H(n)}  \sqcup N_{H(n)}^{+}$
and that this is a triangular decomposition with respect to the 
functional $\Lambda$ given by $\Lambda (\varepsilon_1) = 1$ 
and $\Lambda (\varepsilon_i) = 0$ for $i > 1$. Therefore 
$P_{H(n)}$ is a principal parabolic set of roots and 
its Levi decomposition is $P_{H(n)} =L_{H(n)} \sqcup N_{H(n)}^{+}$.
Using Proposition \ref{appendix2} 
one verifies that $P_{H(n)}$ is a cominuscule
parabolic set of roots.

\begin{theorem} 
\label{prop_hn}
Let $\gg = H(n)$, $n \geq 5$, and $l = \left[ n/2 \right]$. The set 
of all cominuscule parabolic 
subsets of the root system of $H(n)$ forms a single orbit 
under the action of the Weyl group $W_{H(n)_0} \cong W_{\gs \go (n)}$. 
The parabolic subset  $P_{H(n)}$ provides a representative of this orbit. 
Moreover, the following isomorphisms hold for the corresponding Levi subalgebra $\gl$ and nilradical $\gn^+$ 
(considered as an $\gl$-module): $\gl \cong H(n-2) \otimes \bigwedge(1) \oplus \C^2$, ${\gn}^+ \cong \widetilde{H}(n-2) \oplus \C$.
\end{theorem}
\noindent
{\bf Proof.} Let $P$ be a cominuscule parabolic subset of $\Delta$ 
with Levi decomposition $P = L \sqcup N^+$.
If $N^+ \cap \Delta_{\gs \go (n)} =  \emptyset$, then one easily proves that
that $P = L = \Delta$. Assume that $N^+\cap \Delta_{\gs \go (n)} \neq  \emptyset$. We proceed with a case-by-case verification checking for which $i$, $\varepsilon_i$ is in $P$.

\medskip
\noindent {\it Case 1: $n=2l+1$}. It follows from Proposition \ref{algebras} (ii) 
that one can conjugate $P$ by an element of $W_{H(n)_0}$ so that
$L\cap \Delta_{\gs \go (2l+1)}  = L_{\gs \go (2l+1)}$ 
and $N^+ \cap \Delta_{\gs \go (2l+1)}  = N_{\gs \go (2l+1)}^+$. 
This implies that $\varepsilon_1 \in N^+$ and $\varepsilon_i \in L$ for $i>1$. Now one easily obtains that $P = P_{H(2l+1)}$.

\medskip
\noindent {\it Case 2: $n=2l$}.  Following Proposition \ref{algebras} (iv) we proceed 
with three subcases.

\medskip
\noindent {\it Case 2.1: $P\cap \Delta_{\gs \go (2l)} = 
L_{\gs \go (2l)}(1) \sqcup N_{\gs \go (2l)}^+ (1) $}. 
With the aid of Proposition \ref{roots} we easily find 
that $\varepsilon_1 \in N^+$ and $\varepsilon_i \in L$, for all $i>1$. 
Thus $P = P_{H(2l)}$.

\medskip
\noindent {\it Case 2.2: $P\cap \Delta_{\gs \go (2l)} = 
L_{\gs \go (2l)}(l) \cup N_{\gs \go (2l)}^+ (l) $}.  Using again 
Proposition \ref{roots} we verify that $\varepsilon_i \in N^+$ 
for all $i=1,2,...,l$. Since $l>2$, this contradicts to $P$ being cominuscule.

\medskip
\noindent {\it Case 2.3: $P\cap \Delta_{\gs \go (2l)} = 
L_{\gs \go (2l)}(l-1) \cup N_{\gs \go (2l)}^+ (l-1) $}. 
We reach a contradiction in a similar fashion to the previous subcase.

The isomorphisms for the Levi subalgebra and nilradical follow from the explicit description of the root spaces of $H(n)$, see e.g. \cite[Appendix A]{GP2}. More precisely, we have:
\begin{eqnarray*}
\widetilde{\gl} & =& {\mathcal D} \left( \bigwedge(\xi_2, \ldots,\xi_l, \xi_{l+2}, \ldots,
\xi_n ) \oplus \xi_1 \xi_{l+1}  \bigwedge(\xi_2, \ldots ,\xi_l, \xi_{l+2}, 
\ldots,\xi_n ) \right)\\
 \widetilde{\gn}^+ & \cong & \widetilde{H}(\xi_2, \ldots,\xi_l, \xi_{l+2},\ldots,\xi_n ) 
\oplus \C D_{\eta_1},
\end{eqnarray*}
where $ \eta_1:=\frac{1}{\sqrt{2}}(\xi_1 +
\sqrt{-1}\xi_{l+1}),$ and $\widetilde{\gl}$ and $ \widetilde{\gn}^+ $ denote the corresponding to the sets of roots $L_{H(n)}$ and $N_{H(n)}^+$ subalgebras  of $\widetilde{H}(n)$.

\hfill $\square$

\begin{remark}
Theorem \ref{prop_hn} implies that all cominuscule parabolic sets of roots of 
$H(n)$ are principal. One should note though that not every parabolic subset of roots 
of $H(n)$ is principal as shown in \cite[\S 3]{DFG}.
\end{remark}

\begin{remark}
Due to the fact that for every root $\varepsilon_I - \varepsilon_J$ of $H(n)$, 
the intersection $H(n)_i \cap H(n)^{\varepsilon_I - \varepsilon_J}$ is nontrivial 
for more than one index $i$, the subalgebra $H(n)_{-1} \oplus H(n)_0$ is not 
a cominuscule parabolic subalgebra of $H(n)$ according to our definition. 
Some authors studied versions of parabolic subalgebras of simple 
finite dimensional Lie superalgebras, which are not root subalgebras
but at the same time are more restrictive than our definition in different 
respects (see for example \cite{IO} which deals with $\Z$-graded 
Lie superalgebras). It will be interesting
to study cominuscule subalgebras in those frameworks. For instance
$H(n)_{-1} \oplus H(n)_0$ will be an example of such a subalgebra. 
\end{remark}
\section*{Appendix}
\label{appendix}
In this Appendix we prove Proposition \ref{appendix2}.
We need to show that if $\alpha, \beta, \alpha + \beta \in \Delta$, 
then $\left[ \gg^{\alpha}, \gg^{\beta}\right] \neq 0$. 
For all classical Lie superalgebras $\gg$, 
except $\gg = {\bf s} \gp (n), \gp \gs \gq (n)$, 
the proposition follows from Proposition 2.5.5(e) in 
\cite{K}. For the two remaining cases of classical Lie superalgebras and for 
$\gg = W(n), H(n)$, 
it is sufficient to give examples of elements $x_{\alpha} \in \gg^{\alpha}$ and $x_{\beta} \in \gg^{\beta}$ such that 
$[x_{\alpha}, x_{\beta}] \neq 0$. For $\gg = {\bf s} \gp (n), \gg = \gp \gs  \gq (n)$, 
any nonzero $x_\alpha,x_{\beta}$ will work, see \cite{FSS}. 

Let now $\gg = W(n)$. Using the explicit description of $\gg^{\alpha}$ 
provided in \S \ref{W_n}, we provide examples considering three cases 
for $\alpha$ and $\beta$: (i) $\{ \alpha = \varepsilon_{I},  \beta= \varepsilon_{I'}\}$, 
(ii) $\{ \alpha = \varepsilon_{I, j},  \beta= \varepsilon_{I'}, j \notin I\}$, (iii) 
$\{ \alpha = \varepsilon_{I,j},  \beta= \varepsilon_{I',j'}, j\notin I, j' \notin I'\}$.
We only consider case (iii) as the other two are similar. Since $\alpha + \beta \in \Delta$, we may assume 
that $j \in I'$. Then for $I = \{i_1, \ldots ,i_k \}$ 
and $I' = \{ i_1', \ldots ,i_l'\}$, $i_1'=j$, 
we choose $x_{\alpha} = \xi_{i_1} \ldots \xi_{i_k} \frac{\partial}{\partial \xi_j}$ and 
$x_{\beta} = \xi_{i_1'} \ldots \xi_{i_l'} \frac{\partial}{\partial \xi_{j'}}$. 
Then $[x_{\alpha}, x_{\beta}]$ will contain the nonzero homogeneous summand 
$\xi_{i_1} \ldots \xi_{i_k}\xi_{i_2'} \ldots \xi_{i_l'}  
\frac{\partial}{\partial \xi_{j'}}$, 
and hence is nonzero. 

Assume now that $\gg = H(n)$. The description of the root spaces $\gg^{\alpha}$ and the graded root spaces $\gg^{\alpha} \cap \gg_i$ can be found in the Appendix A of \cite{GP2}. The examples of $x_{\alpha}$ 
and $x_{\beta}$ are given according 
to that description on a case-by-case basis. Consider the case 
$n = 2l$, $\alpha = \varepsilon_I - \varepsilon_J$. 
Let $\beta = \varepsilon_{I'} - \varepsilon_{J'}$. 
For convenience we will treat  $I,J,I',J'$ both as ordered tuples 
and as subsets of $[1,n]$. 
Note that $I \cap J = I\cap I' = I' \cap J' = J \cap J' = \emptyset$. 
We first consider the case when either
$I \cap J' \neq \emptyset$ or $I' \cap J \neq \emptyset$.
Let $x_{\alpha} =  D_{\eta_I \eta_{\widehat{J}}}$ 
and $x_{\beta} =  D_{\eta_{I'} \eta_{\widehat{J'}}}$, where 
$\eta_i := \frac{1}{\sqrt{2}} \left( \xi_i + \sqrt{-1}\xi_{i+l}\right)$,  
$\eta_{i+l} := \frac{1}{\sqrt{2}} \left( \xi_i - \sqrt{-1}\xi_{i+l}\right)$, 
$1 \leq i \leq n$, and $\eta_I = \eta_{i_1} \ldots \eta_{i_k}$, 
$\widehat{I} = (i_1+l, \ldots ,i_k + l)$ whenever $I = (i_1, \ldots ,i_k)$, 
$1 \leq i_1< \ldots <i_k \leq l$. Then $[x_{\alpha}, x_{\beta}] = D_{F}$, 
where $F = \{\eta_I \eta_{\widehat{J}},  \eta_{I'} \eta_{\widehat{J'}}\}$. 
Here $\{ f,g\} = (-1)^{\deg f} 
\left( \sum_{i=1}^l \frac{\partial f}{\partial \eta_{i+l}} 
\frac{\partial g}{\partial \eta_{i}} +  
\sum_{i=1}^l   \frac{\partial f}{\partial \eta_{i}} 
\frac{\partial g}{\partial \eta_{i+l}}  \right)$. Using the assumption 
that $I \cap J' \neq \emptyset$ or  $I' \cap J \neq \emptyset$, we verify that 
$[x_{\alpha}, x_{\beta}] \neq 0$. Now consider the case 
$I \cap J'  = I' \cap J = \emptyset$. 
Because $\alpha + \beta \in \Delta$ we also 
have $ I \cap I' = J \cap J'  = \emptyset$. 
Take $k \in J'$. Then $k \notin I, k \notin J$, so we may choose 
$x_{\alpha} = D_{\eta_I \eta_k \eta_{\widehat{J}} \eta_{k+l}}$ 
and again $x_{\beta} =  D_{\eta_{I'} \eta_{\widehat{J'}}}$. We have that 
$\{\eta_I \eta_k\eta_{\widehat{J}} \eta_{k+l}, 
\eta_{I'} \eta_{\widehat{J'}}\} \neq 0$ and hence 
$[x_{\alpha}, x_{\beta}] \neq 0$. The case $n= 2l +1$ 
is treated analogously.
This completes the proof of Proposition \ref{appendix2}. 
\hfill $\square$


\begin{thebibliography}{99}

\bibitem[BL]{BL} S. Billey, V. Lakshmibai, {\em{Singular loci of Schubert varieties}}, 
Prog. Math., 182. Birkh\"auser, Boston, MA, 2000.

\bibitem[Bo]{Bo} N. Bourbaki, {\em{\'El\'ements de math\'ematique. 
Groupes et alg\`ebres de Lie, Ch. IV -- VI}}, Herman, Paris 1968, 288 pp.

\bibitem[BG]{BGY} K. A. Brown, K.R. Goodearl, M. Yakimov, {\em{Poisson structures 
on affine spaces and flag varieties. I. Matrix affine Poisson space}}, 
Adv. Math. {\bf{206}} (2006), 567--629.

\bibitem[DFG]{DFG} I. Dimitrov, V. Futorny, D. Grantcharov, 
{\em{Parabolic sets of roots}}, Contemp. Math. {\bf 499} (2009), 61--73.

\bibitem[FSS]{FSS} L. Frappat, A. Sciarrino, P. Sorba, 
{\em{Dictionary on Lie Algebras and Superalgebras}}, 
Academic Press, San Diego, CA, 2000.

\bibitem[GY]{GY} K. R. Goodearl, M. Yakimov, {\em{Poisson structures on 
affine spaces and flag varieties. II}}, Trans. Amer. Math. Soc. 
{\bf{361}} (2009), 5753--5780.  

\bibitem[GP]{GP} D. Grantcharov, A. Pianzola, 
{\em{Automorphisms and twisted loop algebras 
of finite dimensional simple Lie superalgebras}}, 
Int. Math. Res. Not. {\bf{73}} (2004), 3937--3962.

\bibitem[GP2]{GP2} D. Grantcharov, A. Pianzola, 
{\em{Automorphisms of toroidal Lie superalgebras}}, 
J. Algebra {\bf{319}} (2008),  4230--4248. 

\bibitem[H]{H} S. Helgason, {\em Differential Geometry, Lie Groups, and Symmetric Spaces}, Academic Press, New York, 1978.

\bibitem[IO]{IO} I. Ivanova, A. Onishchik, {\em{Parabolic subalgebras 
and gradings of reductive Lie superalgebras}}, (Russian), 
Sovrem. Mat. Fundam. Napravl. {\bf{20}} (2006), 5--68.

\bibitem[K]{K} V. Kac, {\em{Lie superalgebras}}, 
Adv. Math. {\bf{26}} (1977), 8--96.

\bibitem[LW]{LW} T. Lam, L. Williams, {\em{Total positivity 
for cominuscule Grassmannians}}, New York J. Math. 
{\bf{14}} (2008), 53–99.

\bibitem[P]{P} I. Penkov, {\em{Characters of strongly generic irreducible Lie 
superalgebra representations}}, 
Internat. J. Math.  {\bf{9}}  (1998), 331--366.

\bibitem[PS]{PS} I. Penkov, V. Serganova, 
{\em{Generic irreducible representations 
of finite-dimensional Lie superalgebras}}, 
Internat. J. Math. {\bf{5}} (1994), 389--419.

\bibitem[Po]{Po} A. Postnikov, {\em{Total positivity, Grassmannians, 
and networks}}, preprint arXiv:math/0609764, 79 pp.

\bibitem[PSo]{PSo} K. Purbhoo, F. Sottile, {\em{The recursive nature 
of cominuscule Schubert calculus}}, Adv. Math. {\bf{217}} 
(2008), 1962--2004.

\bibitem[Sch]{Sch} M. Scheunert, {\em{Invariant supersymmetric 
multilinear forms and the Casimir elements of $P$-type 
Lie superalgebras}}, J. Math. Phys. {\bf{28}} (1987), 1180--1191.

\bibitem[S]{S} V. Serganova, {\em{Automorphisms of simple Lie 
superalgebras}}, Math. USSR Izvestiya {\bf{24}} (1985), 539--551.

\bibitem[S2]{S2} V. Serganova, {\em{Kac--Moody superalgebras 
and integrability}}, Developments and trends in infinite-dimensional Lie theory, 169--218,
Progr. Math., 288, Birkh\"auser Boston, Inc., Boston, MA, 2011. 

\bibitem[TY]{TY} H. Thomas, A. Yong, {\em{A combinatorial rule 
for (co)minuscule Schubert calculus}}, 
Adv. Math. {\bf{222}} (2009), 596–-620.
\end{thebibliography}
\end{document}